\documentclass[11pt]{article}
\usepackage{amsfonts}
\usepackage{graphics}
\usepackage{indentfirst}
\usepackage{color}
\usepackage{cite}
\usepackage{latexsym}
\usepackage[paper=a4paper, left=1.6cm, right=1.6cm, top=1.8cm, bottom=1.6cm, headheight=5.5pt, footskip=0.8cm, footnotesep=0.8cm, centering, includefoot]{geometry}
\usepackage{amsmath}
\usepackage{enumitem}
\allowdisplaybreaks
\usepackage{amssymb}
\usepackage[colorlinks, linkcolor=red]{hyperref}
\hypersetup{urlcolor=red, citecolor=red}
\usepackage[dvips]{epsfig}
\usepackage{amscd}
\usepackage{amsthm}
\usepackage{mathrsfs}
\usepackage{verbatim}

\DeclareMathOperator{\divv}{div}
\DeclareMathOperator{\curl}{curl}
\DeclareMathOperator{\loc}{loc}
\allowdisplaybreaks

\begin{document}
\title{Global weak solutions with vacuum for isentropic compressible Navier--Stokes and Navier--Stokes--Poisson equations in the unit disk under Dirichlet boundary conditions
\thanks{Wang's research was partially supported by Postgraduate Research and Innovation Project of Chongqing (No. CYB260103).
Zhong's research was partially supported by National Natural Science Foundation of China (No. 12371227) and
Fundamental Research Funds for the Central Universities (No. SWU--KU24001).
}
}

\author{Shuai Wang,\ Xin Zhong {\thanks{Corresponding author. E-mail addresses: swang238@163.com (S. Wang),
xzhong1014@amss.ac.cn (X. Zhong).}}\date{}\\
\footnotesize
School of Mathematics and Statistics, Southwest University, Chongqing 400715, P. R. China}

\maketitle
\newtheorem{theorem}{Theorem}[section]
\newtheorem{definition}{Definition}[section]
\newtheorem{lemma}{Lemma}[section]
\newtheorem{proposition}{Proposition}[section]
\newtheorem{corollary}{Corollary}[section]
\newtheorem{remark}{Remark}[section]
\renewcommand{\theequation}{\thesection.\arabic{equation}}
\catcode`@=11 \@addtoreset{equation}{section} \catcode`@=12
\maketitle{}

\begin{abstract}
We construct global weak solutions to the isentropic compressible Navier--Stokes and Navier--Stokes--Poisson systems in the unit disk under no-slip boundary conditions, allowing vacuum and unbounded densities for sufficiently small initial energy. Such solutions particularly possess uniform $L^{q_0}$ density bounds for some $q_0>4\gamma$, beyond the natural energy class. Our analysis combines Green function representations with the mean-value property of harmonic functions to isolate a pressure-generated harmonic component of the effective viscous flux and obtain a purely interior potential representation valid for both systems. These results complement Perepelitsa's half-space theory (Arch. Ration. Mech. Anal. 212 (2014), pp. 709--726) and our recent half-plane work (arXiv:2601.11852, 2026), treating bounded-domain geometry and Poisson coupling with nonnegative initial density.
\end{abstract}



%%%%%%%%%%%%%%%%%%%%%%%%%%%%%%%%%%%%%%%%%%%%%%%%%%%%%%%%%%%%%%%%%%%%%%%%%%%%%%%%%%%%%%%%%%%%%%%%%%

\tableofcontents

\section{Introduction}

Let $\Omega\subset\mathbb R^N$ ($N\geq2$) be a bounded domain occupied by a viscous isentropic compressible fluid.
When coupled to a self-consistent field, the fluid motion is described by the Navier--Stokes--Poisson (NSP) system in $\Omega\times(0,\infty)$:
\begin{equation}\label{1.1}
 \begin{cases}
  \rho_t+\divv(\rho\mathbf u)=0,\\
  (\rho\mathbf u)_t+\divv(\rho\mathbf u\otimes\mathbf u)+\nabla P=\mu\Delta\mathbf u+(\mu+\lambda)\nabla\divv\mathbf u+\rho\nabla\phi,\\
  \Delta\phi=\rho-\bar\rho.
 \end{cases}
\end{equation}
Here $\rho=\rho(\mathbf x,t)\geq0$, $\mathbf u=\mathbf u(\mathbf x,t)=(u^1,\ldots,u^N)^\top$, and $\phi=\phi(\mathbf x,t)$
denote the fluid density, velocity, and self-consistent potential, respectively.
The pressure obeys the isentropic law $P=P(\rho)=a\rho^\gamma$ with $a>0$ and $\gamma>1$.
The spatial mean density is given by $\bar\rho=|\Omega|^{-1}\int_\Omega\rho(\mathbf x,t)\mathrm{d}\mathbf{x}$.
The constant viscosity coefficients $\mu$ and $\lambda$ satisfy
\begin{equation*}
  \mu>0,\quad 2\mu+N\lambda\geq0.
\end{equation*}

The system \eqref{1.1} is supplemented with the initial condition
\begin{equation}\label{1.2}
  \rho(\mathbf x,0)=\rho_0(\mathbf x),\quad \rho \mathbf u(\mathbf x,0)=\rho_0\mathbf u_0(\mathbf x), \quad \mathbf x\in\Omega.
\end{equation}
We impose the Dirichlet (no-slip) boundary condition on the velocity and the homogeneous Neumann boundary condition on the potential, given by
\begin{equation}\label{1.3}
  \mathbf u(\mathbf x,t)=\mathbf 0,\quad \partial_{\mathbf n}\phi(\mathbf x,t)=0, \quad (\mathbf x,t)\in\partial\Omega\times(0,\infty),
\end{equation}
where $\mathbf n$ is the unit outward normal to $\partial\Omega$.
The continuity equation \eqref{1.1}$_1$ and the no-slip boundary condition imply that $\bar\rho=|\Omega|^{-1}\int_\Omega\rho_0(\mathbf x)\mathrm{d}\mathbf{x}>0$.
We normalize the potential by requiring $\int_\Omega\phi\mathrm{d}\mathbf{x}=0$ for each $t>0$.

The potential $\phi$ describes the nonlocal interaction induced by the density fluctuation $\rho-\bar\rho$ relative to a uniform neutralizing background.
This interpretation is consistent with Gauss's law in electrostatics, which relates the electric field to the charge distribution \cite{Melcher1981}.
The resulting field acts on the fluid through the force $\rho\nabla\phi$.
For further background on the NSP system and its mathematical formulation, see \cite{Donatelli2003,LMZ2010}. When one neglects the Poisson equation \eqref{1.1}$_3$ and the associated force, \eqref{1.1} reduces to
the isentropic compressible Navier--Stokes (NS) equations
\begin{align}\label{1.4}
 \begin{cases}
  \rho_t+\divv(\rho\mathbf u)=0,\\
  (\rho\mathbf u)_t+\divv(\rho\mathbf u\otimes\mathbf u)+\nabla P=\mu\Delta\mathbf u+(\mu+\lambda)\nabla\divv\mathbf u.
 \end{cases}
\end{align}

\subsection{Background and motivation}\label{sec1.1}

The effective viscous flux plays a central role in the weak-solution theories of Lions--Feireisl \cite{PL98,F04} and Hoff \cite{Hoff95,Hoff95*,Hoff02}.
For the compressible NS and NSP systems, this flux $G:=(2\mu+\lambda)\divv\mathbf{u}-P$ satisfies
\[
  \Delta G=\divv(\rho\dot{\mathbf{u}}-\rho\nabla\phi),
\]
where the Poisson term is absent in the NS case. In the Cauchy setting, this elliptic structure provides a crucial regularity gain for $G$ in Hoff's approach. Under the no-slip condition, however, the absence of a suitable boundary condition for $G$ prevents a direct application of the corresponding elliptic estimates. In particular, the elliptic equation leaves a harmonic component that is not determined by its source term.

The mean-value property of harmonic functions provides a natural link between boundary and volume averages on the unit ball $B$:
\begin{equation}\label{1.5}
  \frac{1}{|\partial B|}\int_{\partial B}h\mathrm{d}S=h(\mathbf{0})=\frac{1}{|B|}\int_B h\mathrm{d}\mathbf{x},\qquad \Delta h=0 \quad \text{in } B.
\end{equation}
From a functional-analytic perspective, the harmonic Bergman projection offers a related viewpoint, as explored in the study of stationary compressible NS boundary value problems \cite{PRS09}; see also \cite{PRS08}. These considerations raise a natural question: can one isolate a harmonic component of the effective viscous flux that captures the boundary effects? We pursue this idea through a Green function representation, exploiting the geometry of the disk to obtain the effective-flux estimates required for global weak solutions of the compressible NS and NSP systems under no-slip boundary conditions.

\subsection{Related work}

Over the past several decades, substantial progress has been made on the global theory of compressible fluid equations. 
Let's review some fundamental developments for the isentropic compressible Navier--Stokes system and its Poisson-coupled counterpart. 
For small perturbations of constant equilibrium states, Matsumura and Nishida \cite{MN80,MN83} established the global well-posedness of classical solutions to the three-dimensional Cauchy and initial-boundary value problems. 
Using Littlewood--Paley analysis, Danchin \cite{Da00} subsequently developed a global well-posedness theory in critical Besov spaces, with initial velocity in $\dot{B}^{N/2-1}_{2,1}(\mathbb{R}^N)$ for $N\geq 2$. 
This framework was further extended by Charve--Danchin \cite{CD10} and independently by Chen--Miao--Zhang \cite{CMZ10} to $L^p$-based critical Besov spaces, allowing velocity data in $\dot{B}^{N/p-1}_{p,1}(\mathbb{R}^N)$ for suitable $p\in[2,2N)$ (see also \cite{H11}). 
For the NSP system, the nonlocal Poisson coupling introduces additional low-frequency effects, which were addressed through hybrid Besov spaces to establish global strong solutions in dimensions $N\geq3$ \cite{HL09}. 
Accommodating vacuum, however, requires analytical techniques beyond these perturbative frameworks.

In the presence of vacuum, the vanishing density complicates the hyperbolic--parabolic structure of the compressible NS system,
making the construction of global solutions substantially more delicate.
Lions \cite{PL98} showed the global existence of finite-energy weak solutions for $\gamma\geq\frac{3N}{N+2}$ in dimensions $N=2,3$,
and Feireisl--Novotn\'y--Petzeltov\'a \cite{FNP01,F04} subsequently generalized the range of adiabatic exponent to $\gamma>\frac{N}{2}$ through refined compactness arguments with the aid of oscillation defect measures.
The effective viscous flux enters these compactness arguments through its weak continuity, which is instrumental in establishing strong convergence of the density.
The resulting weak-solution framework accommodates vacuum, with the density controlled at the natural energy level $L^\infty_tL^\gamma_{\mathbf x}$.
Further developments have addressed more general constitutive laws and special geometric settings.
Bresch and Jabin \cite{BJ18} investigated global weak solutions for non-monotone pressure laws and anisotropic viscous stress tensors,
while Jiang and Zhang \cite{JZ01,JZ03} proved global existence for every $\gamma>1$ under spherical symmetry or axisymmetry.
For the compressible NSP system, the Poisson coupling introduces additional difficulties in the energy estimates and compactness arguments.
Feireisl and collaborators \cite{DFPS04} constructed global weak solutions for self-gravitating flows with possibly non-monotone pressure laws,
and global weak solutions for two-dimensional flows with logarithmically growing pressures were subsequently obtained in \cite{ZT07}.
Despite these advances, uniqueness in the general finite-energy weak-solution class remains an open problem for both the NS and NSP systems.

The pursuit of uniqueness naturally motivates the study of weak solutions with enhanced regularity. 
In this direction, Hoff \cite{Hoff95,Hoff95*,Hoff02} developed a theory of intermediate weak solutions lying between finite-energy weak solutions and strong solutions,
allowing the density to be bounded but discontinuous under suitable assumptions on the initial data.
A key feature of his theory is the interplay between hyperbolic density transport and the elliptic regularity of the effective viscous flux,
which together yield an upper bound for the density along particle trajectories.
This structure also underlies the propagation and persistence of density discontinuities \cite{Hoff02,Hoff08}. Under additional assumptions, Hoff \cite{Hoff06} further proved uniqueness and continuous dependence for such weak solutions. Beyond the Cauchy setting, global weak solutions with vacuum were also constructed in the half-space under Navier-slip boundary conditions \cite{Hoff05}. In parallel, elliptic regularity of the effective viscous flux has played an important role in the \textit{a priori} estimates for global strong solutions with vacuum in the 3D and 2D Cauchy problems \cite{LX19,HLX12} and in initial-boundary value problems under slip boundary conditions \cite{D12,CL23}. A common advantage of suitable slip boundary conditions is that they allow a Neumann-type elliptic treatment of the effective viscous flux, yielding improved spatial regularity.
Since the Poisson force leaves the principal elliptic part of the Lam\'e system unchanged, related methods have also been applied to the NSP system to establish intermediate weak solutions to the Cauchy problem \cite{S2020} and global classical solutions in 3D bounded domains under slip boundary conditions \cite{CHS24}. In contrast, the physically natural no-slip condition does not provide comparable boundary information for the effective viscous flux, calling for a more delicate treatment of boundary effects.

The boundary singularities studied by Hoff and Perepelitsa \cite{HP09,HP12} further illustrate this difficulty. In the half-space, Perepelitsa \cite{P14} extended Hoff's intermediate weak-solution theory to the no-slip condition for small perturbations of a strictly positive equilibrium, assuming a sufficiently small viscosity ratio $\mu/(\mu+\lambda)$ and the boundary-adapted H\"older condition
\begin{equation}\label{1.6}
  [\rho_0]_\alpha^*=\operatorname*{ess\,sup}_{\mathbf{x},\mathbf{y}\in\mathbb{R}^N_+}\frac{|\rho_0(\mathbf{x})-\rho_0(\mathbf{y})|}{|\mathbf{x}-\mathbf{y}^*|^\alpha},
\end{equation}
where $\mathbf{y}^*$ is the reflection of $\mathbf{y}$ across the boundary.
This condition controls the singular integral associated with the reflected Dirichlet Green function
\begin{equation}\label{1.7}
  \mathcal G(\mathbf{x},\mathbf{y})=\Gamma(\mathbf{x}-\mathbf{y})-\Gamma(\mathbf{x}-\mathbf{y}^*),\qquad
  |\nabla^2\Gamma(\mathbf{x}-\mathbf{y}^*)|\lesssim|\mathbf{x}-\mathbf{y}^*|^{-N}.
\end{equation}
The approach was generalized to 3D smooth bounded domains \cite{P18} through Green function estimates and a boundary integral equation.
Since the reflected singular integral does not admit a general $L^\infty$ endpoint estimate, it is natural to seek finite-exponent $L^p$ bounds for the density with $p>\gamma$. 
Such bounds place the pressure in $L^{p/\gamma}$ with an exponent greater than one, 
allowing the singular-integral estimates to be used without imposing additional boundary H\"older regularity.
Along these lines, Huang--Xin--Yan \cite{HXY24} obtained global weak solutions in a 2D disk under a small-energy assumption, using the explicit Dirichlet Green function to reduce the effective-flux boundary remainder to a multiple of $\int_{\partial\Omega}G\mathrm{d}S$. More recently, the current authors \cite{WZ26} constructed global weak solutions in the half-plane with interior and far-field vacuum, exploiting the flat-boundary cancellation $(\mathbf{x}G)\cdot\mathbf{n}=0$ on $\{x_N=0\}$ to derive a Green function decomposition within an $L^p$ density framework.

These results highlight the role of domain geometry in controlling the effective viscous flux under no-slip boundary conditions.
Although the flat-boundary cancellation exploited in \cite{WZ26} is available in the half-space, it does not extend directly to the unit disk.
The key idea of the present work is to exploit a different cancellation mechanism arising from the harmonic structure of the disk.
More precisely, we isolate the pressure-generated harmonic component of the effective viscous flux and use the mean-value property to eliminate its contribution to the boundary-flux integral.
This yields a purely interior potential representation of the effective viscous flux, without the boundary integral arising in \cite{HXY24}, together with the corresponding $L^p$ estimates for both the NS and NSP systems.
Under the assumptions of our main theorems, these estimates enable us to establish global weak solutions allowing vacuum and possessing higher integrability beyond the energy level.

\subsection{Notation and main results}

We fix some notation used throughout the paper.

\begin{itemize}

\item Unless otherwise stated, all function spaces are defined over $\Omega$, with the same notation used for scalar- and vector-valued functions.
For $1\le p\le\infty$ and $k\ge0$, we write
\[
  L^p=L^p(\Omega),\quad W^{k,p}=W^{k,p}(\Omega),\quad H^k=W^{k,2},\quad H_0^1=\overline{C_c^\infty(\Omega)}^{\,H^1}.
\]
The norm in a function space $X$ is denoted by $\|\cdot\|_X$.

\item We use
\[
  \int f\mathrm{d}\mathbf{x}:=\int_\Omega f\mathrm{d}\mathbf{x},\quad f_i:=\partial_i f,\quad
  \bar f:=\frac{1}{|\Omega|}\int_\Omega f\mathrm{d}\mathbf{x},
\]
and denote the material derivative by $\dot f:=f_t+(\mathbf u\cdot\nabla) f$.

\item The effective viscous flux and vorticity are defined by
\begin{equation}\label{1.8}
  F:=(2\mu+\lambda)\divv\mathbf u-(P-\bar P), \quad
  \curl\mathbf{u}:=-\nabla^\bot\cdot\mathbf{u}=-(\partial_2u^1-\partial_1u^2).
\end{equation}
For matrices $A,B\in\mathbb R^{2\times2}$, we write $A:B:=\operatorname{tr}(AB^\top)$.

\item The symbol $C$ denotes a generic positive constant that may vary from line to line,
with its dependence indicated when necessary.
We write $A\lesssim B$ if $A\le CB$.

\end{itemize}

For the NSP system, we define the relative potential energy and the initial total energy by
\begin{equation}\label{1.9}
  G(\rho\,|\,\bar\rho):=\rho\int_{\bar\rho}^{\rho}\frac{P(\xi)-P(\bar\rho)}{\xi^2}\mathrm{d}\xi,\qquad
  E_0:=\int_\Omega\bigg(\frac12\rho_0|\mathbf u_0|^2+G(\rho_0\,|\,\bar\rho)+\frac{1}{2}|\nabla\phi_0|^2\bigg)\mathrm{d}\mathbf{x}.
\end{equation}
Here $\phi_0$ solves $\Delta\phi_0=\rho_0-\bar\rho$ in $\Omega$,
with $\partial_{\mathbf n}\phi_0=0$ on $\partial\Omega$ and $\int_\Omega\phi_0\mathrm{d}\mathbf{x}=0$.
For $\gamma>1$ and $\bar\rho>0$, we have
\begin{equation}\label{1.10}
  G(\rho\,|\,\bar\rho)\ge C
 \begin{cases}
  |\rho-\bar\rho|^2, & |\rho-\bar\rho|\le1,\\ |\rho-\bar\rho|^\gamma, & |\rho-\bar\rho|>1,
 \end{cases}
\end{equation}
where $C>0$ depends only on $a$, $\gamma$, and $\bar\rho$.

The notion of weak solution adopted here is based on the distributional formulations in \cite{P14,S2020}.

\begin{definition}\label{d1.1}
A triplet $(\rho,\mathbf{u},\phi)$ is called a weak solution of \eqref{1.1}--\eqref{1.3} if
\[
  (\rho,\rho\mathbf{u},\nabla\phi)\in C\big([0,\infty);H^{-1}(\Omega)\big),\quad \mathbf{u}\in L^2\big(0,\infty;H_0^1(\Omega)\big),
\]
and
\[
  \phi\in L^2_{\loc}\big([0,\infty);H^1(\Omega)\big),\quad \int_\Omega\phi\mathrm{d}\mathbf{x}=0 \quad\text{for a.e. }t>0.
\]
Moreover, for any $0\le t_1<t_2<\infty$ and any test functions $(\varphi,\boldsymbol\psi)\in C^1\big(\overline{\Omega}\times[t_1,t_2]\big)$
with $\boldsymbol{\psi}=\mathbf 0$ on $\partial\Omega\times[t_1,t_2]$, the following identities hold:
\footnote{We use the Einstein summation convention throughout the paper.}
\begin{align*}
  \left.\int_\Omega\rho\varphi\mathrm{d}\mathbf{x}\right|_{t_1}^{t_2}
  &=\int_{t_1}^{t_2}\int_\Omega(\rho\varphi_t+\rho\mathbf{u}\cdot\nabla\varphi)\mathrm{d}\mathbf{x}\mathrm{d}t,
  \\
  \left.\int_\Omega\rho\mathbf{u}\cdot\boldsymbol{\psi}\mathrm{d}\mathbf{x}\right|_{t_1}^{t_2}
  &=\int_{t_1}^{t_2}\int_\Omega\big(\rho\mathbf{u}\cdot\boldsymbol{\psi}_t+(\rho\mathbf{u}\otimes\mathbf{u}):\nabla\boldsymbol{\psi}
  +P\divv\boldsymbol{\psi}+\rho\nabla\phi\cdot\boldsymbol{\psi}\big)\mathrm{d}\mathbf{x}\mathrm{d}t
  \\ &\quad
  -\int_{t_1}^{t_2}\int_\Omega\big(\mu\nabla\mathbf{u}:\nabla\boldsymbol{\psi}+(\mu+\lambda)\divv\mathbf{u}\divv\boldsymbol{\psi}\big)\mathrm{d}\mathbf{x}\mathrm{d}t,
  \\
  -\int_{t_1}^{t_2}\int_\Omega\nabla\phi\cdot\nabla\varphi\mathrm{d}\mathbf{x}\mathrm{d}t
  &=\int_{t_1}^{t_2}\int_\Omega(\rho-\bar\rho)\varphi\mathrm{d}\mathbf{x}\mathrm{d}t.
\end{align*}
\end{definition}

Our main result is the global existence of weak solutions to the compressible NSP system.

\begin{theorem}\label{t1.1}
Let $\Omega$ be the unit disk. Fix $\gamma\in(1,2)$, $q_0>4\gamma$, $\widehat\rho>\bar\rho+1$, and $M>1$.
There exists $\lambda_*>0$, depending only on $\gamma$, $q_0$, and $\mu$, such that the following holds for every $\lambda\geq\lambda_*$.
Assume that the initial data satisfy
\begin{equation}\label{1.11}
  \rho_0\geq0\quad\text{a.e. in }\Omega,\quad \|\rho_0\|_{L^{q_0}}\leq\widehat\rho,\quad \mathbf u_0\in H_0^1,\quad \|\nabla\mathbf u_0\|_{L^2}\leq M.
\end{equation}
Set
\begin{equation}\label{1.12}
  \vartheta:=\min\left\{\frac{q_0-2\gamma}{8(q_0-\gamma)}+\frac14,\,\frac{\gamma(q_0-4\gamma)}{4(q_0-\gamma)}\right\}.
\end{equation}
There exists a constant $\varepsilon>0$, depending only on $\widehat\rho$, $\bar\rho$, $a$, $\gamma$, $\mu$, $\lambda$, and $q_0$, such that if
\begin{equation}\label{1.13}
  E_0\leq\varepsilon M^{-6/\vartheta},
\end{equation}
then the initial-boundary value problem \eqref{1.1}--\eqref{1.3} admits a global weak solution $(\rho,\mathbf u,\phi)$ in the sense of Definition~\ref{d1.1} satisfying,
for any $T>0$,
\[
  \rho\geq 0\quad\text{a.e. in }\Omega\times(0,T),\quad\rho\in L^\infty(0,\infty;L^{q_0})\cap C([0,T];L^q) \quad(1\leq q<q_0),
\]
and
\begin{align*}
  &\sup_{t\geq0}\int_\Omega\big(\rho|\mathbf u|^2+G(\rho\,|\,\bar\rho)+|\nabla\phi|^2+\sigma|\nabla\mathbf u|^2+\sigma^3\rho|\dot{\mathbf u}|^2\big)\mathrm{d}\mathbf{x}
  \notag\\ &\quad
  +\int_0^\infty\int_\Omega\big(|\nabla\mathbf u|^2+\sigma\rho|\dot{\mathbf u}|^2+\sigma^3|\nabla\dot{\mathbf u}|^2\big)\mathrm{d}\mathbf{x}\mathrm{d}t
  \leq C E_0^\vartheta,
\end{align*}
where $\sigma(t)=\min\{1,t\}$.
\end{theorem}

As a direct application of Theorem \ref{t1.1}, we obtain global weak solutions of the compressible NS system \eqref{1.4} with the initial condition \eqref{1.2} and the no-slip boundary condition in \eqref{1.3}.
Weak solutions are understood as in Definition~\ref{d1.1}, with the Poisson-related requirements omitted.
The corresponding initial energy is
\[
  \widehat E_0:=\int_\Omega\bigg(\frac12\rho_0|\mathbf u_0|^2+G(\rho_0\,|\,\bar\rho)\bigg)\mathrm{d}\mathbf{x}.
\]

\begin{theorem}\label{t1.2}
Let $\Omega, \gamma, q_0, \widehat\rho, M, \vartheta$, and $\lambda_*$ be as in Theorem~\ref{t1.1}, and assume that the initial data satisfy \eqref{1.11}.
There exists a constant $\widehat\varepsilon>0$, depending only on $\widehat\rho, \bar\rho, a, \gamma, \mu, \lambda$, and $q_0$, such that if
\begin{equation}\label{1.14}
  \widehat E_0\leq\widehat\varepsilon M^{-6/\vartheta},
\end{equation}
then the compressible NS system \eqref{1.4} with the initial condition \eqref{1.2} and the no-slip boundary condition in \eqref{1.3} admits a global weak solution $(\rho,\mathbf u)$ verifying, for any $T>0$,
\begin{gather*}
  \rho\geq 0\quad\text{a.e. in }\Omega\times(0,T),\quad\rho\in L^\infty(0,\infty;L^{q_0})\cap C([0,T];L^q) \quad(1\leq q<q_0),\\
  \sup_{t\geq0}\int_\Omega\big(\rho|\mathbf u|^2+G(\rho\,|\,\bar\rho)+\sigma|\nabla\mathbf u|^2+\sigma^3\rho|\dot{\mathbf u}|^2\big)\mathrm{d}\mathbf{x}
  +\int_0^\infty\int_\Omega\big(|\nabla\mathbf u|^2+\sigma\rho|\dot{\mathbf u}|^2+\sigma^3|\nabla\dot{\mathbf u}|^2\big)\mathrm{d}\mathbf{x}\mathrm{d}t
  \leq C \widehat{E}_0^\vartheta.
\end{gather*}
\end{theorem}

We make several remarks on the above results.

\begin{remark}
It should be noted that the weak solutions constructed in Theorems \ref{t1.1} and \ref{t1.2} lie in an intermediate density-integrability regime
between the $L^\gamma$ energy framework of Lions--Feireisl and Hoff's theory \cite{Hoff95,Hoff95*,Hoff02,Hoff05}.
In particular, the uniform $L^{q_0}$ density bounds require neither strict positivity nor an $L^\infty$ upper bound,
thus accommodating both vacuum and unbounded densities.
\end{remark}

\begin{remark}
Compared with Perepelitsa's half-space theory \cite{P14}, the present results allow vacuum without imposing additional boundary H\"older regularity.
Theorem \ref{t1.2} complements the disk result in \cite{HXY24} by relaxing the density-integrability requirement from $q_0=12\gamma/(\gamma-1)$ to any $q_0>4\gamma$.
Unlike our earlier half-plane approach \cite{WZ26}, which relies on a flat-boundary cancellation, the present work addresses the distinct geometry of the unit disk and extends the analysis to the Poisson-coupled system.
\end{remark}

\begin{remark}
A key ingredient common to both systems is the Dirichlet Lam\'e decomposition
\[
  \mathbf{u}=\mathbf{u}_r+\mathbf{u}_s,\quad F=(2\mu+\lambda)\divv\mathbf{u}_r+F_s,\quad F_s:=(2\mu+\lambda)\divv\mathbf{u}_s-(P-\bar P),
\]
where $F_s$ is harmonic with zero mean.
Using the Dirichlet Green function and the mean-value property of harmonic functions in the unit disk, we derive a purely interior representation of $F_s$; see Lemma~\ref{l2.5}.
Since this component is generated solely by the pressure, the resulting representation can be applied to both systems.
We refer to Subsection \ref{sec1.4} for the underlying strategy and Subsection \ref{sec2.2} for the detailed derivation and estimates.
\end{remark}

\begin{remark}
The small-energy conditions \eqref{1.13} and \eqref{1.14} allow large initial velocity gradients, provided the initial energy is correspondingly small.
The restriction $\lambda\geq\lambda_*$, equivalent to a sufficiently small viscosity ratio $\mu/(\mu+\lambda)$, is needed to close the dissipative $L^{q_0}$ density estimates
based on the continuity equation and the effective viscous flux; similar restrictions appear in \cite{P14,P18,HXY24}.
As observed in \cite{P18}, under such conditions, ``the viscous damping of sound waves dominates the growth of amplitude of oscillations due to the reflection of waves from the boundary''.
It is worth noting that the Green function approaches in both the half-space and the unit disk rely on their specific geometries, and extending the present approach to more general domains or more complex fluid models would be an interesting direction for further investigation.
\end{remark}

\subsection{Strategy of the proof}\label{sec1.4}

We outline the main ideas and difficulties in the proof of Theorem \ref{t1.1}. The argument begins with local smooth solutions corresponding to smooth approximations of the initial data with strictly positive densities.
Our aim is to establish \textit{a priori} estimates independent of the positive lower density bound,
which allow global continuation of these solutions and ultimately passage to the vacuum limit by compactness.

A central task is to propagate higher integrability of the density.
While the basic energy inequality provides the natural $L^\gamma$ control, the uniform $L^{q_0}$ bound on $\rho$
supplies an additional integrability needed to handle the pressure and nonlinear terms.
Using the continuity equation \eqref{1.1}$_1$ and the no-slip boundary condition, we obtain
\begin{equation}\label{1.15}
  \frac{\mathrm{d}}{\mathrm{d}t}\int\rho^{q_0}\mathrm{d}\mathbf{x}=-(q_0-1)\int\rho^{q_0}\divv\mathbf{u}\mathrm{d}\mathbf{x}=-\frac{q_0-1}{2\mu+\lambda}\int\rho^{q_0}(F+P-\bar P)\mathrm{d}\mathbf{x}.
\end{equation}
The pressure term has a favorable sign and provides the damping mechanism for the density estimate.
The main difficulty is therefore to control the effective viscous flux $F$, particularly the harmonic contribution arising from the no-slip boundary.

To address this difficulty, we use the Dirichlet Lam\'e decomposition (see \eqref{2.5}--\eqref{2.8})
\[
  F=(2\mu+\lambda)\divv\mathbf{u}_r+F_s,\quad \Delta F_s=0,\quad \int_\Omega F_s\mathrm{d}\mathbf{x}=0.
\]
The forcing-driven component $\mathbf{u}_r$ can be estimated by standard elliptic estimates, whereas the pressure-generated harmonic component $F_s$ has no prescribed boundary data.
To obtain a representation of the latter, we combine the identity $\Delta(\mathbf{x}F_s)=2\nabla F_s$ with the Green function representation.
This leads to the boundary potential in \eqref{2.14}, involving the kernel
\[
  (y_i-x_i)\partial_{x_i}\partial_{\mathbf n_{\mathbf y}}G(\mathbf{x},\mathbf{y}).
\]
The special geometry of the unit disk yields the identity
\[
  (y_i-x_i)\partial_{x_i}\partial_{\mathbf n_{\mathbf y}}G(\mathbf{x},\mathbf{y})=\partial_{\mathbf n_{\mathbf y}}G(\mathbf{x},\mathbf{y})+\frac{1}{2\pi},
  \qquad \mathbf{y}\in\partial\Omega.
\]
By the Poisson representation, the boundary integral in \eqref{2.14} becomes $F_s(\mathbf{x})+F_s(\mathbf{0})$, as shown in \eqref{2.16}.
The term $F_s(\mathbf{x})$ cancels with the corresponding term in the representation, while the mean-value property and the zero-mean condition imply $F_s(\mathbf{0})=0$.
We thus obtain the purely interior representation \eqref{2.9}.
Its kernel is the real part of the Bergman kernel of the unit disk, and the associated $L^p$ boundedness yields \eqref{2.10} for $1<p<\infty$.
Since the representation of $F_s$ involves only the pressure, it has the same form for both the NS and NSP systems.

The displacement weight in this representation is closely related to the classical Lichtenstein construction for Dirichlet Lam\'e systems \cite[pp.~174--177]{MMP1995},
where the factor $\mathbf{y}-\mathbf{x}$ is used to weaken the leading boundary singularity.
Related techniques appear in the half-space through reflected Green functions \cite{P14,WZ26} and in smooth bounded domains through weakly singular boundary integral equations \cite{P18}.
For further background on boundary potentials and integral equations, we refer to \cite{M2000,MMP1995}.
On the unit disk, however, the corresponding kernel representation still produces a boundary term proportional to $\int_{\partial\Omega}F\mathrm{d}s$, as encountered in \cite{HXY24}.
Our key observation is to apply this representation only to the harmonic component $F_s$.
The mean-value property, together with its zero-mean structure, then yields $\int_{\partial\Omega}F_s\mathrm{d}s=0$,
eliminating this boundary contribution without requiring an estimate of the boundary integral of the flux $F$.

We next combine these elliptic estimates with Hoff's time-weighted energy method to control $\sqrt{\rho}\dot{\mathbf{u}}$ and $\nabla\dot{\mathbf{u}}$; see Subsection~\ref{sec3.2}.
Although this method is well established for the Cauchy problem and initial-boundary value problem with slip boundary conditions when an $L^\infty$ density bound is available \cite{Hoff95,Hoff95*,Hoff02,Hoff05},
its application here requires three further ingredients:
\begin{itemize}
    \item[(i)] \textit{Effective-flux estimates under no-slip boundary conditions.}
    The interior representation of $F_s$, together with its harmonic structure, allows us to control its contribution directly in the volume integrals arising from the time-weighted energy estimates.
    In particular, the resulting bounds \eqref{3.20} and \eqref{3.21} avoid estimating the boundary-flux integral encountered in \cite{HXY24}.

    \item[(ii)] \textit{Space-time integrability and propagation of the density bound.}
    In the absence of an $L^\infty$ density bound, the nonlinear terms
    \[
      \sigma'\int_\Omega |P-\bar P||\nabla\mathbf{u}|\mathrm{d}\mathbf{x}+\sigma\int_\Omega |P-\bar P||\nabla\mathbf{u}|^2\mathrm{d}\mathbf{x}
    \]
    require refined space-time estimates.
    We obtain these by treating the short-time and long-time contributions separately; see \eqref{3.36}--\eqref{3.37}.
    To propagate the $L^{q_0}$ density bound, we then use the evolution identity \eqref{1.15} in a Zlotnik-type argument based on Lemma~\ref{l2.4}.
    For sufficiently large $\lambda$, the resulting estimates allow the pressure contribution to control the growth induced by the effective viscous flux; see \eqref{3.54}--\eqref{3.56}.
    This yields the uniform $L^{q_0}$ control needed in the bootstrap argument.

    \item[(iii)] \textit{Poisson coupling.}
    Although the principal Lam\'e operator is unchanged, the nonlocal Poisson force generates additional terms in the energy estimates.
    Their control also requires a space-time $L^2$ estimate for $P-\bar P$.
    To obtain this estimate, we use the pressure evolution derived from the continuity equation in \eqref{3.33}--\eqref{3.35} and
    test the momentum equation \eqref{1.1}$_2$ with a suitable function constructed through the Bogovski\u{\i} operator in Lemma \ref{l2.3}.
    This yields the additional pressure bounds \eqref{3.51}.
\end{itemize}

These estimates are coupled through the bootstrap argument formulated in Proposition \ref{p3.1}.
The assumed $L^{q_0}$ control of the density provides the integrability needed for the time-weighted energy estimates, while the resulting bounds improve the density estimate through \eqref{1.15}.
This closes the bootstrap and yields \textit{a priori} bounds independent of the positive lower density bound.
Consequently, the smooth approximate solutions can be continued globally in time, and compactness allows passage to the vacuum limit, proving Theorem \ref{t1.1} for the NSP system.

The remainder of the paper is organized as follows.
Section \ref{sec2} establishes the preliminary estimates and develops the harmonic decomposition of the effective viscous flux.
Section \ref{sec3} sets up the bootstrap framework and establishes the uniform \textit{a priori} estimates needed to close the bootstrap argument,
and we give the proof of Theorem \ref{t1.1} in Section \ref{sec4}.

\section{Preliminaries}\label{sec2}

In this section we collect auxiliary results and establish the elliptic estimates used later.

\subsection{Auxiliary results and inequalities}

We first recall two interpolation inequalities (see \cite[Theorems 5.8 and 5.9, pp. 139--141]{AF2003}).

\begin{lemma}\label{l2.1}
Let $\Omega\subset\mathbb{R}^N$ be a domain satisfying the cone condition, and let $m\ge1$ be an integer and $1\le r<\infty$.
Then the following inequalities hold.
\begin{enumerate}
\renewcommand{\labelenumi}{\textup{(\roman{enumi})}}
\item
If $mr>N$ and $r\le p\le\infty$, or $mr=N$ and $r\le p<\infty$, or $mr<N$ and $r\le p\le\frac{Nr}{N-mr}$,
then every $f\in W^{m,r}(\Omega)$ satisfies
\begin{equation}\label{2.1}
  \|f\|_{L^p(\Omega)}\le C\|f\|_{W^{m,r}(\Omega)}^\theta\|f\|_{L^r(\Omega)}^{1-\theta},\qquad \theta=\frac{N}{m}\left(\frac1r-\frac1p\right),
\end{equation}
where $1/\infty=0$.

\item
If $r>1$, $mr>N$, and $1\le q<\infty$ satisfies either $q\le r$ or $q>r$ with $(m-1)r<N$,
then every $f\in W^{m,r}(\Omega)\cap L^q(\Omega)$ satisfies
\begin{equation}\label{2.2}
  \|f\|_{L^\infty(\Omega)}\le C\|f\|_{W^{m,r}(\Omega)}^\theta\|f\|_{L^q(\Omega)}^{1-\theta},\qquad \theta=\frac{Nr}{Nr+q(mr-N)}.
\end{equation}
\end{enumerate}
Here $C$ depends only on $m$, $N$, the relevant exponents, and the parameters of the cone condition for $\Omega$.
\end{lemma}

Combining Lemma \ref{l2.1} with the Poincar\'e inequality yields the following Gagliardo--Nirenberg inequalities (cf. \cite[Remark 2.1]{LWZ}).

\begin{lemma}\label{l2.2}
Let $\Omega\subset\mathbb{R}^2$ be a bounded Lipschitz domain, and let $2\le p<\infty$, $1<q<\infty$, and $2<s<\infty$.
For any $f\in H^1(\Omega)$ and $g\in L^q(\Omega)$ with $\nabla g\in L^s(\Omega)$, there hold
\begin{align}\label{2.3}
  \|f\|_{L^p(\Omega)}&\le C_1\|f\|_{L^2(\Omega)}^{\frac2p}\|\nabla f\|_{L^2(\Omega)}^{1-\frac2p}+C_2\|f\|_{L^2(\Omega)},\\
  \|g\|_{L^\infty(\Omega)}&\le C_3\|g\|_{L^q(\Omega)}^{\frac{q(s-2)}{2s+q(s-2)}}\|\nabla g\|_{L^s(\Omega)}^{\frac{2s}{2s+q(s-2)}}+C_4\|g\|_{L^2(\Omega)},\label{2.4}
\end{align}
where $C_i\ge0$ ($i=1,2,3,4$) depend only on $p$, $q$, $s$, and $\Omega$.
Moreover, one may take $C_2=0$ if $f$ has zero mean or, when vector-valued, vanishing normal trace on $\partial\Omega$.
The same holds for $C_4$ with $g$ in place of $f$.
\end{lemma}

To derive a space-time $L^2$ estimate for the pressure fluctuation $P-\bar P$, we use the Bogovski\u{\i} operator (see \cite[Lemma 3.17]{NS04}).

\begin{lemma}\label{l2.3}
Let $\Omega\subset\mathbb{R}^N$ ($N\ge2$) be a bounded Lipschitz domain. For each $1<p<\infty$, there exists a linear operator
\[
  \mathcal B:\overline{L^p}(\Omega)\to W_0^{1,p}(\Omega),\qquad
  \overline{L^p}(\Omega):=\left\{f\in L^p(\Omega):\int_\Omega f\mathrm{d}\mathbf{x}=0\right\},
\]
such that, for all $f\in \overline{L^p}(\Omega)$,
\[
  \divv\mathcal B[f]=f\quad\text{in }\Omega,\qquad
  \|\nabla\mathcal B[f]\|_{L^p(\Omega)}\le C_{p,\Omega}\|f\|_{L^p(\Omega)}.
\]
Moreover, if $f=\divv\mathbf g$ for some $\mathbf g\in L^p(\Omega)$ with vanishing normal trace on $\partial\Omega$, then
\[
  \|\mathcal B[f]\|_{L^p(\Omega)}\le C_{p,\Omega}\|\mathbf g\|_{L^p(\Omega)}.
\]
\end{lemma}

To obtain a uniform-in-time $L^{q_0}$ bound for the density, we use a variant of Zlotnik's lemma (cf. \cite[Lemma 1.3 (b)]{Zlo}).\footnote{The cited result is stated for an equality, but its proof also applies to the differential inequality below.}

\begin{lemma}\label{l2.4}
Let $T>0$, $g\in C(\mathbb R)$, and $y,b\in W^{1,1}(0,T)$ satisfy
\[
  y'(t)\le g\big(y(t)\big)+b'(t) \quad \text{for a.e. }t\in(0,T), \qquad y(0)=y^0.
\]
Suppose that $\lim_{\xi\to\infty}g(\xi)=-\infty$ and, for some $N_0,N_1\ge0$,
\[
  b(t_2)-b(t_1)\le N_0+N_1(t_2-t_1) \quad(0\le t_1<t_2\le T).
\]
Then, for any $\xi^*$ such that $g(\xi)\le-N_1$ whenever $\xi\ge\xi^*$,
\[
  y(t)\le\max\{y^0,\xi^*\}+N_0 \quad(0\le t\le T).
\]
\end{lemma}

\subsection{Harmonic decomposition of the effective viscous flux}\label{sec2.2}

To control the effective viscous flux $F$ under the no-slip boundary condition, we exploit the elliptic structure of the Lam\'e system
to decompose it into a component driven by non-pressure forcing and a harmonic component induced by the pressure gradient.
In the unit disk, the latter admits an explicit interior representation, yielding refined $L^p$ estimates.
Since this harmonic structure is unaffected by the Poisson coupling, the same analysis applies to the compressible NS system.

Writing the momentum equation \eqref{1.1}$_2$ as a Dirichlet Lam\'e system, we have
\begin{equation}\label{2.5}
  \mathcal{L}\mathbf{u}:=\mu\Delta\mathbf{u}+(\mu+\lambda)\nabla\divv\mathbf{u}=\mathbf f+\nabla g,\qquad
  \mathbf f:=\rho\dot{\mathbf{u}}-\rho\nabla\phi,\quad g:=P-\bar P.
\end{equation}
For the NS system, the same formulation holds with $\mathbf f=\rho\dot{\mathbf u}$.
We decompose $\mathbf{u}=\mathbf{u}_r+\mathbf{u}_s$, where
\begin{equation}\label{2.6}
 \begin{cases}
  \mathcal{L}\mathbf{u}_r=\mathbf f,& \mathbf{x}\in\Omega,\\
  \mathbf{u}_r=\mathbf 0,& \mathbf{x}\in\partial\Omega,
 \end{cases}
 \qquad
 \begin{cases}
  \mathcal{L}\mathbf{u}_s=\nabla g,& \mathbf{x}\in\Omega,\\
  \mathbf{u}_s=\mathbf 0,& \mathbf{x}\in\partial\Omega.
 \end{cases}
\end{equation}
Accordingly, the effective viscous flux admits the decomposition
\begin{equation}\label{2.7}
  F=(2\mu+\lambda)\divv\mathbf{u}_r+F_s,\qquad F_s:=(2\mu+\lambda)\divv\mathbf{u}_s-g.
\end{equation}
Taking the divergence of the equation for $\mathbf{u}_s$ and using the zero mean of $g$, we obtain
\begin{equation}\label{2.8}
  \Delta F_s=0\quad\text{in }\Omega,\quad\int_\Omega F_s\mathrm{d}\mathbf{x}=0.
\end{equation}

\begin{lemma}[Decomposition of the effective viscous flux in the unit disk]\label{l2.5}
Let $\Omega$ be the unit disk, and let $(\rho,\mathbf u,\phi)$ be a smooth solution to \eqref{1.1}--\eqref{1.3}.\footnote{The same conclusions hold for the compressible NS system with $\mathbf f=\rho\dot{\mathbf u}$, since the equation defining $\mathbf u_s$ is unchanged.}
Then
\begin{equation}\label{2.9}
  F(\mathbf{x})=(2\mu+\lambda)\divv\mathbf{u}_r(\mathbf{x})+F_s(\mathbf{x}),\quad
  F_s(\mathbf{x})=\frac{2\mu}{3\mu+\lambda}\int_\Omega\partial_{x_i}\partial_{y_i}\Gamma^*(\mathbf{x},\mathbf{y})(P-\bar P)(\mathbf{y})\mathrm{d}\mathbf{y},
\end{equation}
where $G(\mathbf{x},\mathbf{y})=\Gamma(\mathbf{x}-\mathbf{y})-\Gamma^*(\mathbf{x},\mathbf{y})$ is the Dirichlet Green function for the Laplacian in $\Omega$, with $\Delta\Gamma=\delta_0$.
Moreover, for every $1<p<\infty$,
\begin{gather}\label{2.10}
  \|F_s\|_{L^p}\le\frac{C_p\mu}{3\mu+\lambda}\|P-\bar P\|_{L^p},
\\ \label{2.11}
  \|\nabla\mathbf{u}_s\|_{L^p}\le C_{p,\mu,\lambda}\|P-\bar P\|_{L^p},\quad
  \|\mathbf{u}_r\|_{W^{2,p}}\le C_{p,\mu,\lambda}\big(\|\rho\dot{\mathbf{u}}\|_{L^p}+\|\rho\nabla\phi\|_{L^p}\big).
\end{gather}
\end{lemma}

\begin{proof}
Recall that $F_s=(2\mu+\lambda)\divv\mathbf{u}_s-g$ is harmonic in $\Omega$.
From $\mathcal L\mathbf{u}_s=\nabla g$, we have
\[
  \mu\Delta\mathbf{u}_s=-\frac{\mu+\lambda}{2\mu+\lambda}\nabla F_s+\frac{\mu}{2\mu+\lambda}\nabla g.
\]
Since $\Delta(\mathbf{x}F_s)=2\nabla F_s$, it follows that
\begin{equation}\label{2.12}
  \Delta\bigg(\mu\mathbf{u}_s+\frac{\mu+\lambda}{2(2\mu+\lambda)}\mathbf{x}F_s\bigg)=\frac{\mu}{2\mu+\lambda}\nabla g.
\end{equation}

Applying Green's formula to \eqref{2.12} and using $\mathbf{u}_s=\mathbf 0$ on $\partial\Omega$, we obtain
\begin{align}\label{2.13}
  \mu\mathbf{u}_s+\frac{\mu+\lambda}{2(2\mu+\lambda)}\mathbf{x}F_s=
  -\frac{\mu}{2\mu+\lambda}\int_\Omega\nabla_{\mathbf y}G(\mathbf{x},\mathbf{y})g(\mathbf{y})\mathrm{d}\mathbf{y}
  +\frac{\mu+\lambda}{2(2\mu+\lambda)}\int_{\partial\Omega}\partial_{\mathbf n_{\mathbf y}}G(\mathbf{x},\mathbf{y})\mathbf{y}F_s(\mathbf{y})\mathrm{d}S_{\mathbf y}.
\end{align}
Since $G=\Gamma-\Gamma^*$ and
$ -\partial_{x_i}\partial_{y_i}\Gamma(\mathbf{x}-\mathbf{y})=\delta_{\mathbf{x}}(\mathbf{y}),
$
we get
\[
  -\int_\Omega\partial_{x_i}\partial_{y_i}G(\mathbf{x},\mathbf{y})g(\mathbf{y})\mathrm{d}\mathbf{y}
  =g(\mathbf{x})+\int_\Omega\partial_{x_i}\partial_{y_i}\Gamma^*(\mathbf{x},\mathbf{y})g(\mathbf{y})\mathrm{d}\mathbf{y}.
\]
Taking the divergence of \eqref{2.13} and using
\[
  \divv(\mathbf{x}F_s)=2F_s+\mathbf{x}\cdot\nabla F_s
\]
and the Poisson representation of $F_s$, we arrive at
\begin{equation}\label{2.14}
 \boxed{
  F_s(\mathbf{x})=\frac{\mu}{2\mu+\lambda}\int_\Omega\partial_{x_i}\partial_{y_i}\Gamma^*(\mathbf{x},\mathbf{y})g(\mathbf{y})\mathrm{d}\mathbf{y}
  +\frac{\mu+\lambda}{2(2\mu+\lambda)}\int_{\partial\Omega}(y_i-x_i)\partial_{x_i}\partial_{\mathbf n_{\mathbf y}}G(\mathbf{x},\mathbf{y})F_s(\mathbf{y})\mathrm{d}S_{\mathbf y}.
 }
\end{equation}

For the unit disk, the fundamental solution and its harmonic correction are given by
\begin{equation}\label{2.15}
  \Gamma(\mathbf{x}-\mathbf{y})=\frac{1}{2\pi}\log|\mathbf{x}-\mathbf{y}|,\qquad
  \Gamma^*(\mathbf{x},\mathbf{y})=\frac{1}{4\pi}\log\big(1-2\mathbf{x}\cdot\mathbf{y}+|\mathbf{x}|^2|\mathbf{y}|^2\big).
\end{equation}
The corresponding Poisson kernel is
\[
  \partial_{\mathbf n_{\mathbf y}}G(\mathbf{x},\mathbf{y})=\frac{1-|\mathbf{x}|^2}{2\pi|\mathbf{x}-\mathbf{y}|^2},\qquad
  \mathbf{x}\in\Omega,\quad\mathbf{y}\in\partial\Omega.
\]
A direct calculation gives
\[
  (y_i-x_i)\partial_{x_i}\partial_{\mathbf n_{\mathbf y}}G(\mathbf{x},\mathbf{y})=\partial_{\mathbf n_{\mathbf y}}G(\mathbf{x},\mathbf{y})+\frac{1}{2\pi}.
\]
Since $F_s$ is harmonic, its Poisson representation and mean value property yield
\[
  \int_{\partial\Omega}(y_i-x_i)\partial_{x_i}\partial_{\mathbf n_{\mathbf y}}G(\mathbf{x},\mathbf{y})F_s(\mathbf{y})\mathrm{d}S_{\mathbf y}
  =F_s(\mathbf{x})+\frac{1}{2\pi}\int_{\partial\Omega}F_s(\mathbf{y})\mathrm{d}S_{\mathbf y}=F_s(\mathbf{x})+F_s(\mathbf 0).
\]
Substituting this identity into \eqref{2.14} leads to
\begin{equation}\label{2.16}
 \boxed{
  F_s(\mathbf{x})=\frac{2\mu}{3\mu+\lambda}\int_\Omega\partial_{x_i}\partial_{y_i}\Gamma^*(\mathbf{x},\mathbf{y})g(\mathbf{y})\mathrm{d}\mathbf{y}+\frac{\mu+\lambda}{3\mu+\lambda}F_s(\mathbf 0).
 }
\end{equation}
While the mean value property and the zero mean of $F_s$ imply
\begin{equation}\label{2.17}
  F_s(\mathbf 0)=\frac{1}{|\Omega|}\int_\Omega F_s\mathrm{d}\mathbf{x}=0.
\end{equation}
Thus \eqref{2.16} reduces to \eqref{2.9}.

To prove \eqref{2.10}, write $z=x_1+\mathrm{i}x_2$ and $\zeta=y_1+\mathrm{i}y_2$.
Differentiating \eqref{2.15}, we obtain
\[
  \partial_{x_i}\partial_{y_i}\Gamma^*(\mathbf{x},\mathbf{y})=-\frac{1}{\pi}\operatorname{Re}\frac{1}{(1-z\overline{\zeta})^2}.
\]
Since $\pi^{-1}(1-z\overline{\zeta})^{-2}$ is the Bergman kernel of the unit disk, the $L^p$ boundedness of the Bergman projection for $1<p<\infty$ (see \cite{DS2004}) yields
\[
  \left\|\int_\Omega\partial_{x_i}\partial_{y_i}\Gamma^*(\cdot,\mathbf{y})g(\mathbf{y})\mathrm{d}\mathbf{y}\right\|_{L^p}\le C_p\|g\|_{L^p},
\]
which combined with \eqref{2.9} implies \eqref{2.10}.

Finally, the standard $W^{1,p}$ and $W^{2,p}$ estimates for the Dirichlet Lam\'e system, applied respectively to $\mathbf{u}_s$ and $\mathbf{u}_r$, give
\[
  \|\nabla\mathbf{u}_s\|_{L^p}\le C_{p,\mu,\lambda}\|g\|_{L^p},\qquad \|\mathbf{u}_r\|_{W^{2,p}}\le C_{p,\mu,\lambda}\|\mathbf f\|_{L^p}.
\]
These lead to \eqref{2.11}.
\end{proof}

\section{\textit{A priori} estimates}\label{sec3}

Following the approximation procedure in \cite{MN80,S2020}, we consider regularized initial data $(\rho_0^\epsilon,\mathbf u_0^\epsilon)$ satisfying
\begin{equation}\label{3.1}
  \rho_0^\epsilon\in H^2,\qquad \mathbf u_0^\epsilon\in H^2\cap H_0^1,\qquad \inf_\Omega\rho_0^\epsilon>0.
\end{equation}
Let $(\rho^\epsilon,\mathbf u^\epsilon,\phi^\epsilon)$ be a corresponding local strong solution of \eqref{1.1}--\eqref{1.3} on $[0,T_\epsilon]$, satisfying
\begin{equation}\label{3.2}
  (\rho^\epsilon,\mathbf u^\epsilon)\in C\big([0,T_\epsilon];H^2\big),\qquad \phi^\epsilon\in C\big([0,T_\epsilon];H^4\big),\qquad \inf_{\Omega\times[0,T_\epsilon]}\rho^\epsilon>0.
\end{equation}

To carry out the global bootstrap and continuation argument, we establish \textit{a priori} bounds depending on the initial data only through the quantities specified in Theorem \ref{t1.1},
with constants independent of $\epsilon$, the positive lower bound of the density, and the time horizon.
For simplicity, we suppress the superscript $\epsilon$ on the solution and work on $[0,T]$ for any $0<T\le T_\epsilon$. Recalling that $\sigma(t)=\min\{1,t\}$, we define
\begin{align}\label{3.3}
  &\mathcal A_1(T):= \sup_{0\leq t\leq T}\big(\sigma\|\nabla \mathbf{u}\|_{L^2}^2\big)+\int_0^T\sigma\|\sqrt{\rho} \dot{\mathbf{u}}\|_{L^2}^2\mathrm{d}t,\\
  &\mathcal A_2(T):= \sup_{0\leq t\leq T}\big(\sigma^3\|\sqrt{\rho} \dot{\mathbf{u}}\|_{L^2}^2\big)+\int_0^T\sigma^3\|\nabla \dot{\mathbf{u}}\|_{L^2}^2\mathrm{d}t,\label{3.4}\\[0.4em]
  &\mathcal A_3(T):= \sup_{0\leq t\leq T}\|\nabla \mathbf{u}\|_{L^2}^2.\label{3.5}
\end{align}

\subsection{Bootstrap framework}

\begin{proposition}\label{p3.1}
Let the parameters be as in Theorem \ref{t1.1}, with $\lambda\geq\lambda_*$.
There exist constants $\widehat C> 1$ and $\varepsilon>0$, depending only on $\widehat\rho$, $\bar\rho$, $a$, $\gamma$, $\mu$, $\lambda$, and $q_0$,
such that if the strong solution $(\rho,\mathbf u,\phi)$ introduced above exists on $[0,T]$ with initial data satisfying \eqref{1.11} and obeys
\begin{align}\label{3.6}
  \sup_{0\leq t\leq T}\|\rho(t)\|_{L^{q_0}}\leq2\widehat\rho,\quad \mathcal A_1(T)+\mathcal A_2(T)\leq2E_0^{\vartheta},\quad
  \mathcal A_3\big(\sigma(T)\big)\leq3\widehat C M^2,
\end{align}
then these bounds improve to
\begin{align}\label{3.7}
  \sup_{0\leq t\leq T}\|\rho(t)\|_{L^{q_0}}\leq\frac74\widehat\rho,\quad \mathcal A_1(T)+\mathcal A_2(T)\leq E_0^{\vartheta},\quad
  \mathcal A_3\big(\sigma(T)\big)\leq2\widehat C M^2,
\end{align}
provided that $E_0\leq\varepsilon M^{-6/\vartheta}$.
\end{proposition}

The proof of Proposition \ref{p3.1} relies on the \textit{a priori} estimates established below.

\subsection{Energy estimates and bootstrap improvement}\label{sec3.2}

We begin with the basic energy inequality.

\begin{lemma}\label{l3.1}
For any $T>0$, it holds that
\begin{align}\label{3.8}
  \sup_{0\le t\le T}\int\bigg(\frac12\rho|\mathbf{u}|^2+G(\rho\,|\,\bar\rho)+\frac{1}{2}|\nabla\phi|^2\bigg)\mathrm{d}\mathbf{x}
  +\int_0^T\big(\mu\|\nabla\mathbf{u}\|_{L^2}^2+(\mu+\lambda)\|\divv\mathbf{u}\|_{L^2}^2\big)\mathrm{d}t\le 2E_0.
\end{align}
\end{lemma}

\begin{proof}
Multiplying the momentum equation $\eqref{1.1}_2$ by $\mathbf{u}$, integrating over $\Omega$,
and using the continuity equation $\eqref{1.1}_1$ and the definition of $G(\rho\,|\,\bar\rho)$ in \eqref{1.9}, we obtain
\begin{align}\label{3.9}
  \frac{\mathrm{d}}{\mathrm{d}t}\int\bigg(\frac12\rho|\mathbf{u}|^2+G(\rho\,|\,\bar\rho)\bigg)\mathrm{d}\mathbf{x}
  +\mu\|\nabla\mathbf{u}\|_{L^2}^2+(\mu+\lambda)\|\divv\mathbf{u}\|_{L^2}^2=\int\rho\mathbf{u}\cdot\nabla\phi\mathrm{d}\mathbf{x}.
\end{align}
By \eqref{1.1}$_1$ and the Poisson equation \eqref{1.1}$_3$, integration by parts yields
\[
  \int\rho\mathbf{u}\cdot\nabla\phi\mathrm{d}\mathbf{x}=\int\phi\rho_t\mathrm{d}\mathbf{x}=-\frac12\frac{\mathrm{d}}{\mathrm{d}t}\|\nabla\phi\|_{L^2}^2.
\]
Combining this identity with \eqref{3.9} and integrating in time yields \eqref{3.8}.
\end{proof}

We next establish a bound for $\mathcal A_3(\sigma(T))$.

\begin{lemma}\label{l3.2}
Under the bootstrap assumption \eqref{3.6}, there exists a constant $\varepsilon_1>0$, depending only on $\widehat\rho$, $\bar\rho$, $a$, $\gamma$, $\mu$, $\lambda$, and $q_0$,
such that if $E_0\leq\varepsilon_1 M^{-4/(q_0-3)}$, then
\begin{equation}\label{3.10}
  \mathcal A_3\big(\sigma(T)\big)+\int_0^{\sigma(T)}\|\sqrt{\rho}\dot{\mathbf u}\|_{L^2}^2\mathrm{d}t\leq 2\widehat C M^2.
\end{equation}
\end{lemma}

\begin{proof}
We first derive some bounds for the density fluctuation and the Poisson potential.
By \eqref{1.10}, the energy inequality \eqref{3.8}, and H\"older's inequality, one has
\begin{align}\label{3.11}
  \|\rho-\bar\rho\|_{L^\gamma}^\gamma
  &=\int_{\{|\rho-\bar\rho|\leq1\}}|\rho-\bar\rho|^\gamma\mathrm{d}\mathbf{x}+\int_{\{|\rho-\bar\rho|>1\}}|\rho-\bar\rho|^\gamma\mathrm{d}\mathbf{x}
  \notag\\ &\leq
  |\Omega|^{1-\frac{\gamma}{2}}\bigg(\int_{\{|\rho-\bar\rho|\le1\}}|\rho-\bar\rho|^2\mathrm{d}\mathbf{x}\bigg)^{\frac{\gamma}{2}}
  +\int_{\{|\rho-\bar\rho|>1\}}|\rho-\bar\rho|^\gamma\mathrm{d}\mathbf{x}
  \leq C\Big(E_0^{\frac{\gamma}{2}}+E_0\Big).
\end{align}
Interpolating between \eqref{3.11} and the density bound in \eqref{3.6}, we obtain that, if $E_0\le1$,
\begin{align}\label{3.12}
  \|\rho-\bar\rho\|_{L^{q\gamma}}^\gamma
  \le\|\rho-\bar\rho\|_{L^\gamma}^{\frac{\gamma(q_0-q\gamma)}{q(q_0-\gamma)}}
  \|\rho-\bar\rho\|_{L^{q_0}}^{\frac{\gamma q_0(q-1)}{q(q_0-\gamma)}}
  \le CE_0^{\frac{\gamma(q_0-q\gamma)}{2q(q_0-\gamma)}},
  \quad 1\le q\le4,
\end{align}
Then the pressure law gives
\begin{align}\label{3.13}
  \|P-\bar P\|_{L^q}^q&=\int\bigg|P(\rho)-\frac{1}{|\Omega|}\int_\Omega P(\rho)\mathrm{d}\mathbf{y}\bigg|^q\mathrm{d}\mathbf{x}
  \notag\\[0.4em] &\leq
  C\|P-P(\bar\rho)\|_{L^q}^q \leq C\left(\|\rho-\bar\rho\|_{L^{q\gamma}}^{q\gamma}+\|\rho-\bar\rho\|_{L^q}^q\right),\quad 2\le q\le4.
\end{align}
Moreover, the Neumann elliptic estimate and the Sobolev embedding imply
\begin{equation}\label{3.14}
  \|\nabla\phi\|_{L^\infty}\le C\|\phi\|_{W^{2,2\gamma}}\le C\|\rho-\bar\rho\|_{L^{2\gamma}}\le CE_0^{\frac{q_0-2\gamma}{4(q_0-\gamma)}}.
\end{equation}

We turn to the velocity gradient.
The momentum equation \eqref{1.1}$_2$ can be written as
\begin{equation}\label{3.15}
  \rho\dot{\mathbf u}+\nabla(P-\bar P)=\mu\Delta\mathbf u+(\mu+\lambda)\nabla\divv\mathbf u+\rho\nabla\phi.
\end{equation}
Testing \eqref{3.15} against $\mathbf u_t$ and integrating by parts, we find
\begin{align}\label{3.16}
  &\frac12\frac{\mathrm{d}}{\mathrm{d}t}\left(\mu\|\nabla\mathbf u\|_{L^2}^2+(\mu+\lambda)\|\divv\mathbf u\|_{L^2}^2\right)
  +\int\rho|\dot{\mathbf u}|^2\mathrm{d}\mathbf{x}
  \notag\\ &
  =\int\rho\mathbf u_t\cdot\nabla\phi\mathrm{d}\mathbf{x}+\int(P-\bar P)\divv\mathbf u_t\mathrm{d}\mathbf{x}
  +\int\rho\dot{\mathbf u}\cdot(\mathbf u\cdot\nabla)\mathbf u\mathrm{d}\mathbf{x}.
\end{align}
For the Poisson term, using $\mathbf u_t=\dot{\mathbf u}-(\mathbf u\cdot\nabla)\mathbf u$ and \eqref{3.14}, we obtain
\begin{align}\label{3.17}
  \left|\int\rho\mathbf u_t\cdot\nabla\phi\mathrm{d}\mathbf{x}\right|
  \le C\|\nabla\phi\|_{L^\infty}\big(\|\sqrt\rho\dot{\mathbf u}\|_{L^2}+\|\nabla\mathbf u\|_{L^2}^2\big)
  \le\frac18\|\sqrt\rho\dot{\mathbf u}\|_{L^2}^2+C\big(\|\nabla\mathbf u\|_{L^2}^2+1\big).
\end{align}
For the pressure term, Jensen's inequality gives $(\bar P/a)^{1/\gamma}\ge\bar\rho$, and hence
one gets from the triangle inequality and the monotonicity of the pressure law that
\begin{equation}\label{3.18}
  \|\rho-\bar\rho\|_{L^q}^q\le C_q\left\|\rho-\left(\frac{\bar P}{a}\right)^{1/\gamma}\right\|_{L^q}^q
  \le C_{q,a,\bar\rho}\|P-\bar P\|_{L^q}^q, \quad 1\le q<\infty,
\end{equation}
which together with the interpolation inequality \eqref{2.1} and \eqref{2.11} gives
\begin{align}\label{3.19}
  \|\nabla\mathbf u_r\|_{L^4}^4\le C\|\nabla\mathbf u_r\|_{L^2}\|\nabla^2\mathbf u_r\|_{L^{3/2}}^3
  &\le C\left(\|\nabla\mathbf u\|_{L^2}+\|\nabla\mathbf u_s\|_{L^2}\right)
  \left(\|\sqrt\rho\|_{L^6}^3\|\sqrt\rho\dot{\mathbf u}\|_{L^2}^3+\|\rho\|_{L^3}^3\|\nabla\phi\|_{L^3}^3\right)
  \notag\\[0.2em] &\le
  C\left(\|\nabla\mathbf u\|_{L^2}+\|P-\bar P\|_{L^2}\right)\left(\|\sqrt\rho\dot{\mathbf u}\|_{L^2}^3+\|P-\bar P\|_{L^2}^3\right).
\end{align}
Thus, by the continuity equation \eqref{1.1}$_1$ and integrating by parts,
\begin{align}\label{3.20}
  &\int (P-\bar{P})\divv\mathbf{u}_t\mathrm{d}\mathbf{x}
  \notag\\ &=
  \frac{\mathrm{d}}{\mathrm{d}t}\int(P-\bar{P})\divv\mathbf{u}\mathrm{d}\mathbf{x}
  +\int P'(\rho)\big(\mathbf{u}\cdot\nabla\rho+\rho\divv\mathbf{u}\big)\divv\mathbf{u}\mathrm{d}\mathbf{x}
  \notag\\ &=
  \frac{\mathrm{d}}{\mathrm{d}t}\int(P-\bar{P})\divv\mathbf{u}\mathrm{d}\mathbf{x}
  +\int\mathbf{u}\cdot\nabla\big(P-\bar{P}\big)\divv\mathbf{u}\mathrm{d}\mathbf{x}+\int P'(\rho)\rho(\divv\mathbf{u})^2\mathrm{d}\mathbf{x}
  \notag\\ &=
  \frac{\mathrm{d}}{\mathrm{d}t}\int(P-\bar{P})\divv\mathbf{u}\mathrm{d}\mathbf{x}
  +\int\big[P'(\rho)\rho-(P-\bar{P})\big]\big(\divv(\mathbf{u}_r+\mathbf{u}_s)\big)^2\mathrm{d}\mathbf{x}
  -\int \big(P-\bar{P}\big)\mathbf{u}\cdot\nabla\divv(\mathbf{u}_r+\mathbf{u}_s)\mathrm{d}\mathbf{x}
  \notag\\ &\leq
  \frac{\mathrm{d}}{\mathrm{d}t}\int(P-\bar{P})\divv\mathbf{u}\mathrm{d}\mathbf{x}
  +2\int\big[P'(\rho)\rho-(P-\bar{P})\big](\divv\mathbf{u}_r)^2\mathrm{d}\mathbf{x}
  -\int \big(P-\bar{P}\big)\mathbf{u}\cdot\nabla\divv\mathbf{u}_r\mathrm{d}\mathbf{x}
  \notag\\ &\quad
  +\frac{C}{(2\mu+\lambda)^2}\int\big[P'(\rho)\rho-(P-\bar{P})\big]\big[F_s^2+\big(P-\bar P\big)^2\big]\mathrm{d}\mathbf{x}
  +\frac{1}{2(2\mu+\lambda)}\int\big(P-\bar{P}\big)^2\divv\mathbf{u}\mathrm{d}\mathbf{x}
  \notag\\ &\quad
  -\frac{1}{2\mu+\lambda}\int\big(P-\bar{P}\big)\mathbf{u}\cdot\nabla F_s\mathrm{d}\mathbf{x}
  \notag\\ &\leq
  \frac{\mathrm{d}}{\mathrm{d}t}\int(P-\bar{P})\divv\mathbf{u}\mathrm{d}\mathbf{x}
  +C\|P-\bar P\|_{L^4}\|\nabla\mathbf{u}_r\|_{L^{4}}\|\nabla\mathbf{u}_r\|_{L^{2}}
  +C\big(\|\nabla\mathbf{u}_r\|_{L^2}^2+\|\divv\mathbf{u}\|_{L^2}^2\big)
  \notag\\[0.4em] &\quad
  +C\|P-\bar P\|_{L^4}\|\mathbf{u}\|_{L^\frac{4\gamma}{\gamma-1}}\|\nabla^2\mathbf{u}_r\|_{L^\frac{4\gamma}{1+2\gamma}}
  +C\big(\|P-\bar P\|_{L^2}^2+\|P-\bar P\|_{L^3}^3+\|P-\bar P\|_{L^4}^4\big)
  \notag\\ &\quad
  -\frac{1}{2\mu+\lambda}\int\big(P-\bar{P}\big)\mathbf{u}\cdot\nabla F_s\mathrm{d}\mathbf{x}
  \notag\\ &\leq
  \frac{\mathrm{d}}{\mathrm{d}t}\int(P-\bar{P})\divv\mathbf{u}\mathrm{d}\mathbf{x}+\frac18\|\sqrt{\rho}\dot{\mathbf{u}}\|_{L^2}^2
  +C\big(\|\nabla\mathbf{u}\|_{L^2}^2+1\big)-\frac{1}{2\mu+\lambda}\int\big(P-\bar{P}\big)\mathbf{u}\cdot\nabla F_s\mathrm{d}\mathbf{x},
\end{align}
provided that $E_0\le1$, where we have used
\[
  P'(\rho)\rho-(P-\bar{P})=(\gamma-1)(P-\bar{P})+\gamma \bar{P}\leq C\big(P-\bar{P}\big)+C.
\]
To control the last integral in \eqref{3.20}, we invoke the weighted derivative estimate for harmonic functions in \cite[Theorem 5.1]{CKY2001}.
For any harmonic function $h$ in the unit ball $B\subset\mathbb R^N$, $1\le p<\infty$, and integer $m\ge1$, it holds that
\[
  \sum_{|\alpha|<m}|\partial^\alpha h(\mathbf 0)|+\sum_{|\alpha|=m}\big\|(1-|\mathbf{x}|^2)^m\partial^\alpha h\big\|_{L^p(B)}\le C\|h\|_{L^p(B)},
\]
where $\alpha$ denotes a multi-index.
Since $d(\mathbf{x}):=\operatorname{dist}(\mathbf{x},\partial\Omega)=1-|\mathbf{x}|\le1-|\mathbf{x}|^2$,
the harmonicity of $F_s$ and Lemma \ref{l2.5} imply
\[
  \|d\nabla F_s\|_{L^4}\le C\|F_s\|_{L^4}\le C\|P-\bar P\|_{L^4},
\]
which together with Hardy's inequality yields
\[
  \left|\int(P-\bar P)\mathbf u\cdot\nabla F_s\mathrm{d}\mathbf{x}\right|
  \le C\|P-\bar P\|_{L^4}\|d\nabla F_s\|_{L^4}\left\|\frac{\mathbf u}{d}\right\|_{L^2}\le C\|P-\bar P\|_{L^4}^2\|\nabla\mathbf u\|_{L^2}.
\]
Substituting this bound into \eqref{3.20}, we conclude that
\begin{align}\label{3.21}
  \int(P-\bar P)\divv\mathbf u_t\mathrm{d}\mathbf{x}
  \le\frac{\mathrm{d}}{\mathrm{d}t}\int(P-\bar P)\divv\mathbf u\mathrm{d}\mathbf{x}+\frac18\|\sqrt\rho\dot{\mathbf u}\|_{L^2}^2
  +C\left(\|\nabla\mathbf u\|_{L^2}^2+1\right).
\end{align}
It remains to estimate the convection term.
By the interpolation inequality \eqref{2.1} and \eqref{2.11},
\begin{align*}
  \|\nabla\mathbf{u}_r\|_{L^{\frac{2q_0}{q_0-2}}}
  &\le C\|\nabla\mathbf{u}_r\|_{L^2}^{\frac{q_0-3}{q_0-1}}\|\nabla^2\mathbf{u}_r\|_{L^{\frac{2q_0}{q_0+1}}}^{\frac{2}{q_0-1}}
  \notag\\ &\le
  C\|\nabla\mathbf{u}_r\|_{L^2}^{\frac{q_0-3}{q_0-1}}\Big(\|\sqrt{\rho}\|_{L^{2q_0}}^{\frac{2}{q_0-1}}\|\sqrt{\rho}\dot{\mathbf{u}}\|_{L^2}^{\frac{2}{q_0-1}}+1\Big)
  \le C\|\nabla\mathbf{u}_r\|_{L^2}^{\frac{q_0-3}{q_0-1}}\Big(\|\sqrt{\rho}\dot{\mathbf{u}}\|_{L^2}^{\frac{2}{q_0-1}}+1\Big).
\end{align*}
Using H\"older's inequality and the decomposition of $\mathbf u$, we then find
\begin{align}\label{3.23}
  \int\rho \dot{\mathbf{u}}\cdot(\mathbf{u}\cdot\nabla)\mathbf{u}\mathrm{d}\mathbf{x}
  &\leq \|\sqrt{\rho}\dot{\mathbf{u}}\|_{L^2}\|\rho^{1/2}\|_{L^{2q_0}}
  \Big(\|\mathbf{u}\|_{L^{2q_0}}\|\nabla\mathbf{u}_r\|_{L^\frac{2q_0}{q_0-2}}+\|\mathbf{u}\|_{L^\frac{4q_0}{q_0-2}}\|\nabla\mathbf{u}_s\|_{L^4}\Big)
  \notag\\ &\leq
  C\big(\|\sqrt{\rho}\dot{\mathbf{u}}\|_{L^2}+1\big)^\frac{q_0+1}{q_0-1}\|\nabla\mathbf{u}\|_{L^2}\big(\|\nabla\mathbf{u}\|_{L^2}+1\big)^{\frac{q_0-3}{q_0-1}}
  +C\|\sqrt{\rho}\dot{\mathbf{u}}\|_{L^2}\|\nabla\mathbf{u}\|_{L^2}\|P-\bar P\|_{L^4}
  \notag\\ &\leq
  \frac14\|\sqrt{\rho}\dot{\mathbf{u}}\|_{L^2}^2+C\Big(\|\nabla\mathbf{u}\|_{L^2}^{\frac{4(q_0-2)}{q_0-3}}+1\Big).
\end{align}

Consequently, substituting \eqref{3.17}, \eqref{3.21}, and \eqref{3.23} into \eqref{3.16}, we arrive at
\begin{equation}\label{3.24}
  \frac12\frac{\mathrm{d}}{\mathrm{d}t}\mathcal E_1(t)+\frac12\|\sqrt\rho\dot{\mathbf u}\|_{L^2}^2
  \le C\bigg(\|\nabla\mathbf u\|_{L^2}^{\frac{4(q_0-2)}{q_0-3}}+1\bigg),
\end{equation}
where
\begin{equation}\label{3.25}
  \mathcal E_1(t):=\mu\|\nabla\mathbf u\|_{L^2}^2+(\mu+\lambda)\|\divv\mathbf u\|_{L^2}^2-2\int(P-\bar P)\divv\mathbf u\mathrm{d}\mathbf{x}
\end{equation}
satisfies
\begin{equation}\label{3.26}
  c_0\|\nabla\mathbf u\|_{L^2}^2-CE_0^{\frac{\gamma(q_0-2\gamma)}{2(q_0-\gamma)}}\le\mathcal E_1(t)
  \le C_0\|\nabla\mathbf u\|_{L^2}^2+CE_0^{\frac{\gamma(q_0-2\gamma)}{2(q_0-\gamma)}},
\end{equation}
due to \eqref{3.12}, \eqref{3.13}, and Young's inequality,
where $c_0,C_0>0$ depend only on the viscosity coefficients.
Integrating \eqref{3.24} over $(0,\sigma(T))$ and combining the result with \eqref{3.26} and \eqref{3.8}, we deduce
\begin{align}\label{3.27}
  \mathcal A_3\big(\sigma(T)\big)+\int_0^{\sigma(T)}\|\sqrt\rho\dot{\mathbf u}\|_{L^2}^2\mathrm{d}t
  &\le C_1M^2+C\int_0^{\sigma(T)}\|\nabla\mathbf u\|_{L^2}^{\frac{4(q_0-2)}{q_0-3}}\mathrm{d}t
  \notag\\ &\le
  C_1M^2+C\big[\mathcal A_3\big(\sigma(T)\big)\big]^{\frac{q_0-1}{q_0-3}}\int_0^{\sigma(T)}\|\nabla\mathbf u\|_{L^2}^2\mathrm{d}t
  \notag\\ &\le
  C_1M^2+C_2E_0\big(3\widehat C M^2\big)^{\frac{q_0-1}{q_0-3}}.
\end{align}
Choosing $\widehat C=\max\left\{2,C_1,27C_2\right\}$ and $\varepsilon_1=\min\left\{1, \widehat{C}^{-\frac{q_0-1}{q_0-3}}\right\}$,
the right-hand side of \eqref{3.27} is bounded by $2\widehat C M^2$ whenever $E_0\le\varepsilon_1M^{-4/(q_0-3)}$.
This proves \eqref{3.10}.
\end{proof}

We now turn to the estimates for $\mathcal A_1(T)$ and $\mathcal A_2(T)$.

\begin{lemma}\label{l3.3}
Suppose that the bootstrap bounds in \eqref{3.6} hold and $\lambda\geq\lambda_1$, where $\lambda_1>0$ depends only on $\gamma$ and $\mu$.
Then there exists $\varepsilon_2>0$, depending only on $\widehat\rho$, $\bar\rho$, $a$, $\gamma$, $\mu$, $\lambda$, and $q_0$, such that
\begin{equation}\label{3.28}
  \mathcal A_1(T)+\mathcal A_2(T)\leq E_0^\vartheta,
\end{equation}
provided that $E_0\leq\varepsilon_2M^{-2}$.
\end{lemma}
\begin{proof}
\textbf{Step 1. Preliminary estimates for $\mathcal A_1(T)$.}
Multiplying \eqref{3.15} by $\sigma\dot{\mathbf u}$ and integrating the resultant over $\Omega$, one gets
\begin{align}\label{3.29}
  \int\sigma\rho|\dot{\mathbf u}|^2\mathrm{d}\mathbf{x}&=-\sigma\int\dot{\mathbf u}\cdot\nabla(P-\bar P)\mathrm{d}\mathbf{x}
  +\sigma\int\big[\mu\Delta\mathbf u\cdot\dot{\mathbf u}+(\mu+\lambda)\dot{\mathbf u}\cdot\nabla\divv\mathbf u\big]\mathrm{d}\mathbf{x}+\sigma\int\rho\dot{\mathbf u}\cdot\nabla\phi\mathrm{d}\mathbf{x}
  \notag\\ &
  =:\mathcal I_1+\mathcal I_2+\sigma\int\rho\dot{\mathbf u}\cdot\nabla\phi\mathrm{d}\mathbf{x}.
\end{align}
Using the continuity equation \eqref{1.1}$_1$, the pressure law $P=a\rho^\gamma$, and integration by parts, we find
\begin{align}\label{3.30}
  \mathcal I_1&=\sigma\int(P-\bar{P})\divv\mathbf{u}_t\mathrm{d}\mathbf{x}-\sigma\int(\mathbf{u}\cdot\nabla)\mathbf{u}\cdot\nabla (P-\bar{P})\mathrm{d}\mathbf{x}
  \notag\\ &=
  \frac{\mathrm{d}}{\mathrm{d}t}\int \sigma (P-\bar{P})\divv\mathbf{u}\mathrm{d}\mathbf{x}-\sigma'\int (P-\bar{P})\divv\mathbf{u}\mathrm{d}\mathbf{x}
  +\sigma\int\big[\mathbf{u}\cdot\nabla(P-\bar{P})+P'(\rho)\rho\divv\mathbf{u}\big]\divv\mathbf{u}\mathrm{d}\mathbf{x}
  \notag\\ &\quad
  -\sigma\int(\mathbf{u}\cdot\nabla)\mathbf{u}\cdot\nabla (P-\bar{P})\mathrm{d}\mathbf{x}
  \notag\\ &=
  \frac{\mathrm{d}}{\mathrm{d}t}\int\sigma (P-\bar{P})\divv\mathbf{u}\mathrm{d}\mathbf{x}
  -\sigma'\int (P-\bar{P})\divv\mathbf{u}\mathrm{d}\mathbf{x}+\sigma\int(P-\bar{P})\partial_iu^j\partial_ju^i\mathrm{d}\mathbf{x}
  \notag\\ &\quad
  +\sigma\int \big[P'(\rho)\rho-(P-\bar{P})\big](\divv\mathbf{u})^2\mathrm{d}\mathbf{x}
  \notag\\ &\leq
  \frac{\mathrm{d}}{\mathrm{d}t}\int\sigma(P-\bar{P})\divv\mathbf{u}\mathrm{d}\mathbf{x}+C\|\nabla\mathbf{u}\|^2_{L^2}
  +C\sigma'\int|P-\bar{P}||\nabla\mathbf{u}|\mathrm{d}\mathbf{x}
  +C\sigma\int|P-\bar{P}||\nabla\mathbf{u}|^2\mathrm{d}\mathbf{x}.
\end{align}
Similarly, integration by parts gives
\begin{align}\label{3.31}
  \mathcal I_2&=-\frac{1}{2}\frac{\mathrm{d}}{\mathrm{d}t}\big[\mu\sigma\|\nabla\mathbf{u}\|^2_{L^2}+(\mu+\lambda)\sigma\|\divv\mathbf{u}\|^2_{L^2}\big]
  +\frac{\sigma'}{2}\big[\mu\|\nabla\mathbf{u}\|^2_{L^2}+(\mu+\lambda)\|\divv\mathbf{u}\|^2_{L^2}\big]
  \notag\\ &\quad
  -\mu\sigma \int \partial_iu^j\partial_i(u^k\partial_ku^j)\mathrm{d}\mathbf{x}-(\mu+\lambda)\sigma \int \partial_iu^i\partial_j(u^k\partial_ku^j)\mathrm{d}\mathbf{x}
  \notag\\ &\leq
  -\frac{1}{2}\frac{\mathrm{d}}{\mathrm{d}t}\big[\mu\sigma\|\nabla\mathbf{u}\|^2_{L^2}+(\mu+\lambda)\sigma\|\divv\mathbf{u}\|^2_{L^2}\big]
  +C\|\nabla\mathbf{u}\|^2_{L^2}+C\sigma\|\nabla\mathbf{u}\|^3_{L^3}.
\end{align}
Combining \eqref{3.29}--\eqref{3.31} and recalling the definition of $\mathcal E_1(t)$ in \eqref{3.25}, we arrive at
\begin{align}\label{3.32}
  \frac12\frac{\mathrm{d}}{\mathrm{d}t}\big(\sigma\mathcal E_1(t)\big)+\sigma\|\sqrt\rho\dot{\mathbf u}\|_{L^2}^2
  &\le\sigma\int\rho\dot{\mathbf u}\cdot\nabla\phi\mathrm{d}\mathbf{x}
  +C\|\nabla\mathbf u\|_{L^2}^2+C\sigma\|\nabla\mathbf u\|_{L^3}^3
  \notag\\ &\quad
  +C\sigma'\int|P-\bar P||\nabla\mathbf u|\mathrm{d}\mathbf{x}+C\sigma\int|P-\bar P||\nabla\mathbf u|^2\mathrm{d}\mathbf{x}.
\end{align}

We next exploit the coupling between the Poisson and continuity equations to estimate the Poisson term.
Differentiating $\eqref{1.1}_3$ in time gives $\Delta\phi_t=-\divv(\rho\mathbf u)$, with $\partial_{\mathbf n}\phi_t=0$ on $\partial\Omega$.
Testing this equation against $\phi_t$ and integrating by parts, we have
\begin{equation}\label{3.33}
  \|\nabla\phi_t\|_{L^2}^2=-\int\rho\mathbf u\cdot\nabla\phi_t\mathrm{d}\mathbf{x}\le\frac12\|\nabla\phi_t\|_{L^2}^2+C\|\nabla\mathbf u\|_{L^2}^2.
\end{equation}
Moreover, the continuity equation $\eqref{1.1}_1$ gives $\rho\dot{\mathbf u}=(\rho\mathbf u)_t+\divv(\rho\mathbf u\otimes\mathbf u)$.
Thus, integration by parts and the no-slip boundary condition yield
\begin{align}\label{3.34}
  \int\rho\dot{\mathbf u}\cdot\nabla\phi\mathrm{d}\mathbf{x}&=\frac{\mathrm{d}}{\mathrm{d}t}\int\rho\mathbf u\cdot\nabla\phi\mathrm{d}\mathbf{x}
  -\int\rho\mathbf u\cdot\nabla\phi_t\mathrm{d}\mathbf{x}-\int(\rho\mathbf u\otimes\mathbf u):\nabla^2\phi\mathrm{d}\mathbf{x}
  \notag\\ &\le
  \frac{\mathrm{d}}{\mathrm{d}t}\int\rho\mathbf u\cdot\nabla\phi\mathrm{d}\mathbf{x}+\|\nabla\phi_t\|_{L^2}^2
  +C\|\rho\|_{L^{q_0}}\|\nabla^2\phi\|_{L^{q_0}}\|\mathbf u\|_{L^{\frac{2q_0}{q_0-2}}}^2
  \notag\\ &\le
  \frac{\mathrm{d}}{\mathrm{d}t}\int\rho\mathbf u\cdot\nabla\phi\mathrm{d}\mathbf{x}+C\|\nabla\mathbf u\|_{L^2}^2.
\end{align}
Integrating \eqref{3.34} multiplied by $\sigma$ over $(0,T)$, we deduce from \eqref{3.8} and \eqref{3.14} that
\begin{align}\label{3.35}
  \int_0^T\sigma\int\rho\dot{\mathbf{u}}\cdot\nabla\phi\mathrm{d}\mathbf{x}\mathrm{d}t
  &\leq CE_0+\bigg|\int_0^{\sigma(T)}\sigma'\int\rho\mathbf{u}\cdot\nabla\phi\mathrm{d}\mathbf{x}\mathrm{d}t\bigg|
  +\sup_{0\leq t\leq T}\bigg|\int\sigma\rho\mathbf{u}\cdot\nabla\phi\mathrm{d}\mathbf{x}\bigg|
  \notag\\ &\leq
  CE_0+C\|\nabla\phi\|_{L^\infty}\bigg(\int_0^{\sigma(T)}\|\nabla\mathbf{u}\|_{L^2}^2\mathrm{d}t\bigg)^\frac12+C\sup_{0\leq t\leq T}\big(\|\nabla\phi\|_{L^\infty}\|\sqrt{\rho}\mathbf u\|_{L^2}\big)
  \notag\\ &\leq
  CE_0+CE_0^{\frac{q_0-2\gamma}{4(q_0-\gamma)}+\frac12}.
\end{align}

The term involving $\sigma'$ in \eqref{3.32} is treated similarly.
Indeed, \eqref{3.8} and the bounds \eqref{3.12}--\eqref{3.13} give
\begin{align}\label{3.36}
  \int_0^T\int\sigma' |P-\bar{P}||\nabla\mathbf{u}|\mathrm{d}\mathbf{x}\mathrm{d}t\leq
  C\|P-\bar{P}\|_{L^2}\bigg(\int_0^{\sigma(T)}\|\nabla\mathbf{u}\|_{L^2}^2\mathrm{d}t\bigg)^\frac12
  \leq CE_0^{\frac{\gamma(q_0-2\gamma)}{4(q_0-\gamma)}+\frac12}.
\end{align}
Integrating \eqref{3.32} over $(0,T)$ and using \eqref{3.35}, \eqref{3.36}, and \eqref{3.19}, one derives
\begin{align}\label{3.37}
  &\sup_{0\leq t\leq T}\big(\sigma\mathcal E_1(t)\big)+\int_0^T\sigma\|\sqrt{\rho}\dot{\mathbf{u}}\|_{L^2}^2\mathrm{d}t
  \notag\\ &\leq
  CE_0^{\frac{q_0-2\gamma}{4(q_0-\gamma)}+\frac12}+C\int_0^T\sigma\|\nabla\mathbf{u}\|^3_{L^3}\mathrm{d}t+C\int_0^T\sigma\int|P-\bar{P}||\nabla\mathbf{u}|^2\mathrm{d}\mathbf{x}\mathrm{d}t
  \notag\\ &\leq
  CE_0^{\frac{q_0-2\gamma}{4(q_0-\gamma)}+\frac12}+C\int_0^{T}\sigma\|\nabla\mathbf{u}\|_{L^2}\big(\|\nabla\mathbf{u}_r\|_{L^4}^2+\|\nabla\mathbf{u}_s\|_{L^4}^2
  +\|\nabla\mathbf{u}_r\|_{L^6}\|P-\bar P\|_{L^3}+\|\nabla\mathbf{u}_s\|_{L^4}\|P-\bar P\|_{L^4}\big)\mathrm{d}t
  \notag\\ &\leq
  CE_0^{\frac{q_0-2\gamma}{4(q_0-\gamma)}+\frac12}+C\int_0^{T}\|\nabla\mathbf{u}\|_{L^2}^2\mathrm{d}t
  +C\sup_{0\leq t\leq T}\big(\sigma\|\nabla \mathbf{u}\|_{L^2}^2\big)\sup_{0\leq t\leq T}\|\nabla \mathbf{u}\|_{L^2}^2\int_0^{T}\|\nabla\mathbf{u}\|_{L^2}^2\mathrm{d}t
  \notag\\ &\quad
  +\frac12\int_0^T\sigma\|\sqrt{\rho}\dot{\mathbf{u}}\|_{L^2}^2\mathrm{d}t+C\int_0^{T}\sigma\|P-\bar P\|_{L^4}^4\mathrm{d}t
  \notag\\ &\leq
  CE_0^{\frac{q_0-2\gamma}{4(q_0-\gamma)}+\frac12}+CE_0+C\mathcal A_1(T)\big[\mathcal A_3\big(\sigma (T)\big)+\mathcal A_1(T)\big]E_0
  +\frac12\int_0^T\sigma\|\sqrt{\rho}\dot{\mathbf{u}}\|_{L^2}^2\mathrm{d}t
  \notag\\ &\quad
  +C\sup_{0\leq t\leq T}\|P-\bar P\|_{L^4}^4\int_0^{\sigma(T)}1\mathrm{d}t+C\int_{\sigma(T)}^{T}\|P-\bar P\|_{L^4}^4\mathrm{d}t
  \notag\\ &\leq
  CE_0^{2\vartheta}+C_3M^2E_0\mathcal A_1(T)+\frac12\int_0^T\sigma\|\sqrt{\rho}\dot{\mathbf{u}}\|_{L^2}^2\mathrm{d}t+C\int_0^{T}\sigma^3\|P-\bar P\|_{L^4}^4\mathrm{d}t,
\end{align}
where $\mathcal E_1(t)$ and $\vartheta$ are defined in \eqref{3.25} and \eqref{1.12}, respectively.
Consequently, if $E_0\leq\varepsilon_{2,1}M^{-2}$ with $\varepsilon_{2,1}:=\min\{\varepsilon_1,(2C_3)^{-1}\}$, then
\begin{equation}\label{3.38}
  \mathcal A_1(T)\leq CE_0^{2\vartheta}+C\int_0^{T}\sigma^3\|P-\bar P\|_{L^4}^4\mathrm{d}t.
\end{equation}

\textbf{Step 2. Preliminary estimates for $\mathcal A_2(T)$.}
For $m\geq1$, applying $\partial_t+\divv(\mathbf{u}\,\cdot\,)$ componentwise to \eqref{3.15},
multiplying the resulting $j$th equation by $\sigma^m\dot{u}^j$, summing over $j$, and integrating over $\Omega$, we obtain
\begin{align}\label{3.39}
  &\frac{1}{2}\frac{\mathrm{d}}{\mathrm{d}t}\int\sigma^m\rho|\dot{\mathbf{u}}|^2\mathrm{d}\mathbf{x}-\frac{m}{2}\sigma^{m-1}\sigma'\int\rho|\dot{\mathbf{u}}|^2\mathrm{d}\mathbf{x}
  \notag\\ &=
  -\sigma^m\int\dot{u}^j\big[\partial_j P_{t}+\divv(\mathbf{u}\partial_{j}P)\big]\mathrm{d}\mathbf{x}
  +\sigma^m\int\dot u^j\big[\partial_t(\rho\partial_j\phi)+\divv(\rho\mathbf{u}\partial_j\phi)\big]\mathrm{d}\mathbf{x}
  \notag\\ &\quad
  +\mu\sigma^m\int\dot{u}^j\big[\Delta u_{t}^j+\divv\big( \mathbf{u}\Delta u^{j}\big)\big]\mathrm{d}\mathbf{x}
  +(\mu+\lambda)\sigma^m\int\dot{u}^j\big[\partial_j\divv\mathbf{u}_t +\divv(\mathbf{u}\partial_j\divv\mathbf{u})\big]\mathrm{d}\mathbf{x}
  \notag\\ &=:\mathcal{J}_1+\mathcal{J}_2+\mathcal{J}_3+\mathcal{J}_4.
\end{align}
By the continuity equation \eqref{1.1}$_1$, the pressure satisfies
\begin{equation}\label{3.40}
  P_t+\divv(P\mathbf{u})+(\gamma-1)P\divv\mathbf{u}=0,
\end{equation}
which, upon taking the spatial average, yields
\begin{equation}\label{3.41}
  (P-\bar P)_t+\mathbf{u}\cdot\nabla(P-\bar P)+\gamma P\divv\mathbf{u}-(\gamma-1)\overline{P\divv\mathbf{u}}=0.
\end{equation}
Using the continuity equation \eqref{1.1}$_1$ and the identity \eqref{3.40}, we take advantage of \eqref{3.14} and \eqref{3.33} to estimate $\mathcal J_1$ and $\mathcal J_2$ as follows
\begin{align}\label{3.42}
  \mathcal{J}_1&=\sigma^m\int P_{t}\divv\dot{\mathbf{u}}\mathrm{d}\mathbf{x}
  -\sigma^m\int\dot{\mathbf{u}}\cdot\nabla\divv(P\mathbf{u})\mathrm{d}\mathbf{x}+\sigma^m\int\dot{u}^j\divv(P\partial_{j}\mathbf{u})\mathrm{d}\mathbf{x}
  \notag\\ &=
  \sigma^m\int \big(P_{t}+\divv(P\mathbf{u})\big)\divv\dot{\mathbf{u}}\mathrm{d}\mathbf{x}
  +\sigma^m\int(\dot{\mathbf{u}}\cdot\nabla)\mathbf{u}\cdot\nabla P\mathrm{d}\mathbf{x}+\sigma^m\int P\dot{\mathbf{u}}\cdot\nabla\divv\mathbf{u}\mathrm{d}\mathbf{x}
  \notag\\ &=
  -\sigma^m\int(\gamma-1)(P-\bar P+\bar P)\divv\mathbf{u}\divv\dot{\mathbf{u}}\mathrm{d}\mathbf{x}-\sigma^m\int (P-\bar P+\bar P)\partial_i\dot{u}^j\partial_ju^i \mathrm{d}\mathbf{x}
  \notag\\ &\leq
  \frac{\mu}{8}\sigma^m\|\nabla\dot{\mathbf{u}}\|_{L^2}^2+C\|\nabla\mathbf{u}\|_{L^2}^2+C\sigma^m\|\nabla\mathbf{u}\|_{L^4}^4+C\sigma^m\|P-\bar P\|_{L^4}^4,
  \\ \label{3.43}
  \mathcal J_2&=\sigma^m\int\rho\dot u^j\big(\partial_j\phi_t+u^k\partial_k\partial_j\phi\big)\mathrm{d}\mathbf{x}
  \notag\\ &\leq
  C\sigma^m\|\rho\|_{L^4}\|\dot{\mathbf{u}}\|_{L^4}\Big(\|\nabla \phi_t\|_{L^2}+\|\mathbf{u}\|_{L^\frac{2\gamma}{\gamma-1}}\|\nabla^2 \phi\|_{L^{2\gamma}}\Big)
  \notag\\ &\leq
  \frac{\mu}{8}\sigma^m\|\nabla\dot{\mathbf{u}}\|_{L^2}^2+C\|\nabla\mathbf{u}\|_{L^2}^2.
\end{align}
For the viscous terms, integration by parts yields
\begin{align}\label{3.44}
  \mathcal{J}_3&=\mu\sigma^m\int\dot{u}^j\big[\Delta\dot{u}^j-\Delta(\mathbf{u}\cdot\nabla u^j)+\divv\big( \mathbf{u}\Delta u^{j}\big)\big]\mathrm{d}\mathbf{x}
  \notag\\ &=
  \mu\sigma^m\int\big[-|\nabla\dot{\mathbf{u}}|^2+\dot{u}^j_i(u^ku_k^j)_i-\dot{u}^j_i(u^ku^j_i)_k-\dot{u}^j(u^k_iu^j_i)_k\big]\mathrm{d}\mathbf{x}
  \notag\\ &=
  \mu\sigma^m\int\big[-|\nabla\dot{\mathbf{u}}|^2+\dot{u}^j_i(u^ku_k^j)_i-\dot{u}^j_i(u^ku^j_i)_k+\dot{u}_k^j(u^k_iu^j_i)\big]\mathrm{d}\mathbf{x}
  \notag\\ &\leq
  -\frac{7\mu}{8}\sigma^m\|\nabla\dot{\mathbf{u}}\|_{L^2}^2+C\sigma^m\|\nabla\mathbf{u}\|_{L^4}^4,
  \\ \label{3.45}
  \mathcal{J}_4&=(\mu+\lambda)\sigma^m\int\dot{u}^j\big[\partial_j\divv\mathbf{u}_t+ \partial_j\divv(\mathbf{u}\divv\mathbf{u})-\divv(\partial_j\mathbf{u}\divv\mathbf{u})\big]\mathrm{d}\mathbf{x}
  \notag\\ &=
  -(\mu+\lambda)\sigma^m\int\divv\dot{\mathbf{u}}\big[\divv\mathbf{u}_t+\divv(\mathbf{u}\divv\mathbf{u})\big]\mathrm{d}\mathbf{x}
  -(\mu+\lambda)\sigma^m\int\dot{u}^j\divv(\partial_j\mathbf{u}\divv\mathbf{u})\mathrm{d}\mathbf{x}
  \notag\\ &=
  -(\mu+\lambda)\sigma^m\int\big[\divv\dot{\mathbf{u}}\big(\divv\dot{\mathbf{u}}-\partial_iu^j\partial_ju^i+(\divv\mathbf{u})^2\big)-\partial_j\mathbf{u}\cdot\nabla\dot{u}^j\divv\mathbf{u}\big]\mathrm{d}\mathbf{x}
  \notag\\&\leq
  -\frac{\mu+\lambda}{2}\sigma^m\|\divv\dot{\mathbf{u}}\|_{L^2}^2+\frac{\mu}{8}\sigma^m\|\nabla\dot{\mathbf{u}}\|_{L^2}^2+C\sigma^m\|\nabla\mathbf{u}\|_{L^4}^4.
\end{align}
Thus, inserting \eqref{3.42}--\eqref{3.45} into \eqref{3.39}, one sees
\begin{equation}\label{3.46}
  \frac{\mathrm{d}}{\mathrm{d}t}\int\sigma^m\rho|\dot{\mathbf{u}}|^2\mathrm{d}\mathbf{x}+\mu\sigma^m\|\nabla\dot{\mathbf{u}}\|_{L^2}^2
  \leq m\sigma^{m-1}\sigma'\|\sqrt{\rho}\dot{\mathbf{u}}\|_{L^2}^2+C\|\nabla\mathbf{u}\|_{L^2}^2+C\sigma^m\|\nabla\mathbf{u}\|_{L^4}^4+C\sigma^m\|P-\bar P\|_{L^4}^4.
\end{equation}
Taking $m=3$ and integrating \eqref{3.46} over $(0,T)$, we derive from \eqref{3.4}, \eqref{3.8}, and \eqref{3.19} that
\begin{align}\label{3.47}
  \mathcal A_2(T)&\leq C\int_0^T\sigma\|\sqrt{\rho}\dot{\mathbf{u}}\|_{L^2}^2\mathrm{d}t+CE_0
  +C\int_0^{T}\sigma^3\big(\|\nabla \mathbf{u}_r\|_{L^4}^4+\|\nabla \mathbf{u}_s\|_{L^4}^4\big)\mathrm{d}t+C\int_0^{T}\sigma^3\|P-\bar P\|_{L^4}^4\mathrm{d}t
  \notag\\ &\leq
  C\mathcal A_1(T)+CE_0+C\sup_{0\leq t\leq T}\big[\big(\sigma^\frac12\|\nabla\mathbf{u}\|_{L^2}\big)
  \big(\sigma^\frac32\|\sqrt{\rho}\dot{\mathbf{u}}\|_{L^2}\big)\big]
  \int_0^T\sigma\|\sqrt{\rho}\dot{\mathbf{u}}\|_{L^2}^2\mathrm{d}t+C\int_0^{T}\sigma^3\|P-\bar P\|_{L^4}^4\mathrm{d}t
  \notag\\ &\leq
  CE_0+C\mathcal A_1(T)+C\mathcal A_1^\frac32(T)\mathcal A_2^\frac12(T)+C\int_0^{T}\sigma^3\|P-\bar P\|_{L^4}^4\mathrm{d}t.
\end{align}

\textbf{Step 3. Closing the estimates for $\mathcal A_1(T)$ and $\mathcal A_2(T)$.}
Combining \eqref{3.38} and \eqref{3.47}, one infers from the bootstrap bounds in \eqref{3.6} that
\begin{align}\label{3.48}
  \mathcal A_1(T)+\mathcal A_2(T)\leq CE_0^{2\vartheta}+C\mathcal A_1(T)+C\int_0^{T}\sigma^3\|P-\bar P\|_{L^4}^4\mathrm{d}t
  \leq CE_0^{2\vartheta}+C\int_0^{T}\sigma^3\|P-\bar P\|_{L^4}^4\mathrm{d}t.
\end{align}
To control the pressure contribution, we multiply \eqref{3.41} by $3\sigma^3(P-\bar P)^2$ and integrate over $\Omega$. Then it follows from \eqref{2.10} and Young's inequality that
\begin{align*}
  \frac{3\gamma-1}{2\mu+\lambda}\sigma^3\|P-\bar{P}\|_{L^4}^4&=-\frac{\mathrm{d}}{\mathrm{d}t}\int\sigma^3(P-\bar{P})^3\mathrm{d}\mathbf{x}
  +3\sigma^2\sigma'\|P-\bar{P}\|_{L^3}^3+3(\gamma-1)\overline{P\divv\mathbf{u}}\int\sigma^3(P-\bar{P})^2\mathrm{d}\mathbf{x}
  \notag\\ &\quad
  -\frac{3\gamma-1}{2\mu+\lambda}\int\sigma^3(P-\bar{P})^3\big[(2\mu+\lambda)\divv\mathbf{u}_r+F_s\big]\mathrm{d}\mathbf{x}
  -3\gamma \bar{P}\int\sigma^3(P-\bar{P})^2\divv\mathbf{u}\mathrm{d}\mathbf{x}
  \notag\\ &\leq
  -\frac{\mathrm{d}}{\mathrm{d}t}\int\sigma^3(P-\bar{P})^3\mathrm{d}\mathbf{x}+C\sigma'\|P-\bar{P}\|_{L^3}^3
  +\frac{3\gamma-1}{4(2\mu+\lambda)}\left(1+\frac{C_\gamma\mu}{3\mu+\lambda}\right)\left\|P-\bar{P}\right\|_{L^4}^4
  \notag\\[0.2em] &\quad
  +C\sigma^3\|\nabla\mathbf{u}_r\|_{L^4}^4+C\|\nabla\mathbf{u}\|_{L^2}^2.
\end{align*}
Choose $\lambda\geq\lambda_1$ so that $C_\gamma\mu/(3\mu+\lambda_1)\leq1$.
Integrating the preceding inequality over $(0,T)$ and using \eqref{3.8}, \eqref{3.12}, \eqref{3.13}, and \eqref{3.19}, we get
\begin{align}\label{3.49}
  \int_0^T\sigma^3\|P-\bar{P}\|_{L^4}^4\mathrm{d}t
  &\leq C\sup_{0\leq t\leq T}\left\|P-\Bar{P}\right\|_{L^3}^3+CE_0 +C\int_0^T\sigma^3\|P-\bar{P}\|_{L^2}^4\mathrm{d}t
   \notag\\ &\quad
  +C\sup_{0\leq t\leq T}\Big[\big(\sigma^\frac12\|\nabla\mathbf{u}\|_{L^2}\big)
  \big(\sigma^\frac32\|\sqrt{\rho}\dot{\mathbf{u}}\|_{L^2}\big)\Big]\int_0^T\sigma\|\sqrt{\rho}\dot{\mathbf{u}}\|_{L^2}^2\mathrm{d}t
  \notag\\ &\leq
  CE_0^{\frac{\gamma(q_0-3\gamma)}{2(q_0-\gamma)}}+CE_0+C\sup_{0\leq t\leq T}\left\|P-\Bar{P}\right\|_{L^2}^2\int_0^T\|P-\bar{P}\|_{L^2}^2\mathrm{d}t
  +C\mathcal A_1^\frac32(T)\mathcal A_2^\frac12(T)
  \notag\\ &\leq
  CE_0^{2\vartheta}+CE_0^{\frac{\gamma(q_0-2\gamma)}{2(q_0-\gamma)}}\int_0^T\|P-\bar{P}\|_{L^2}^2\mathrm{d}t.
\end{align}

To estimate the space-time $L^2$ norm of $P-\bar P$, we use the Bogovski\u{\i} operator from Lemma \ref{l2.3}.
Set
\[
  \mathbf{w}:=\mathcal B\left[P-\bar P+C_*(\rho-\bar\rho)\right],
\]
where $C_*>0$ will be chosen below.
By the monotonicity of $P(\rho)=a\rho^\gamma$ and the zero-mean property of $\rho-\bar\rho$, we have
\[
  \int(P-\bar P)\divv\mathbf{w}\mathrm{d}\mathbf{x}=\|P-\bar P\|_{L^2}^2+C_*\int\big(P-P(\bar\rho)\big)(\rho-\bar\rho)\mathrm{d}\mathbf{x}
  \geq\|P-\bar P\|_{L^2}^2.
\]
Testing the momentum equation \eqref{1.1}$_2$ against $\mathbf w$ and integrating by parts, we find
\begin{align}\label{3.50}
  \|P-\bar P\|_{L^2}^2&\leq\frac{\mathrm{d}}{\mathrm{d}t}\int\rho\mathbf{u}\cdot\mathbf{w}\mathrm{d}\mathbf{x}
  -\int\rho\mathbf{u}\cdot\mathbf{w}_t\mathrm{d}\mathbf{x}-\int(\rho\mathbf{u}\otimes\mathbf{u}):\nabla\mathbf{w}\mathrm{d}\mathbf{x}
  \notag\\ &\quad
  +\mu\int\nabla\mathbf{u}:\nabla\mathbf{w}\mathrm{d}\mathbf{x}+(\mu+\lambda)\int\divv\mathbf{u}\divv\mathbf{w}\mathrm{d}\mathbf{x}
  -\int\rho\nabla\phi\cdot\mathbf{w}\mathrm{d}\mathbf{x}
  \notag\\ &=
  \frac{\mathrm{d}}{\mathrm{d}t}\int\rho\mathbf{u}\cdot\mathbf{w}\mathrm{d}\mathbf{x}-\int\rho\mathbf{u}\cdot\mathbf{w}_t\mathrm{d}\mathbf{x}
  -\int(\rho\mathbf{u}\otimes\mathbf{u}):\nabla\mathbf{w}\mathrm{d}\mathbf{x}+\mu\int\nabla\mathbf{u}:\nabla\mathbf{w}\mathrm{d}\mathbf{x}
  \notag\\ &\quad
  +(\mu+\lambda)\int\divv\mathbf{u}\divv\mathbf{w}\mathrm{d}\mathbf{x}+\bar\rho\int\phi(P-\bar P)\mathrm{d}\mathbf{x}
  -C_*\bar\rho\|\nabla\phi\|_{L^2}^2-\int(\rho-\bar\rho)\nabla\phi\cdot\mathbf{w}\mathrm{d}\mathbf{x}
  \notag\\ &\leq
  \frac{\mathrm{d}}{\mathrm{d}t} \int\rho\mathbf{u}\cdot\mathbf{w}\mathrm{d}\mathbf{x}+C\|\nabla\mathbf{u}\|_{L^2}^2
  +\frac14\|P-\bar P\|_{L^2}^2-(C_*\bar\rho-C)\|\nabla\phi\|_{L^2}^2+C_4E_0^{\frac{\gamma(q_0-2)}{4(q_0-\gamma)}}\|P-\bar P\|_{L^2}^2
  \notag\\ &\leq
  \frac{\mathrm{d}}{\mathrm{d}t}\int\rho\mathbf{u}\cdot\mathbf{w}\mathrm{d}\mathbf{x}+C\|\nabla\mathbf{u}\|_{L^2}^2+\frac12\|P-\bar P\|_{L^2}^2,
\end{align}
provided that $C_*$ is sufficiently large and $E_0\leq\varepsilon_{2,2}:=\min\left\{1,(4C_4)^{-\frac{4(q_0-\gamma)}{\gamma(q_0-2)}}\right\}$.
Indeed, the equation \eqref{3.40} and the continuity equation \eqref{1.1}$_1$ give
\[
  \mathbf{w}_t=-\mathcal B\left[\divv\big((P+C_*\rho)\mathbf{u}\big)\right]-(\gamma-1)\mathcal B\left[P\divv\mathbf{u}-\overline{P\divv\mathbf{u}}\,\right].
\]
The Bogovski\u{\i} estimates in Lemma \ref{l2.3} together with \eqref{3.14}, \eqref{3.18}, and the Sobolev embedding yield
\begin{gather*}
  \|\mathbf{w}_t\|_{L^2}\leq C\big(\|(P+C_*\rho)\mathbf{u}\|_{L^2}+\|P\divv\mathbf{u}\|_{L^{4/3}}\big)\leq C\|\nabla\mathbf{u}\|_{L^2},
  \\
  \bigg|\int(\rho-\bar\rho)\nabla\phi\cdot\mathbf{w}\mathrm{d}\mathbf{x}\bigg|
  \leq C\|\rho-\bar\rho\|_{L^2}\|P-\bar P\|_{L^2}^2\leq CE_0^{\frac{\gamma(q_0-2)}{4(q_0-\gamma)}}\|P-\bar P\|_{L^2}^2.
\end{gather*}
Integrating \eqref{3.50} over $(0,T)$ and using the energy inequality \eqref{3.8}, we obtain
\begin{equation}\label{3.51}
  \int_0^T \|P-\bar P\|_{L^2}^2\mathrm{d}t\leq
  C\sup_{0\leq t\leq T}\|\sqrt{\rho}\mathbf u\|_{L^2}+C\int_0^T \|\nabla\mathbf{u}\|_{L^2}^2\mathrm{d}t\leq CE_0^\frac12,
\end{equation}
which combined with \eqref{3.48} and \eqref{3.49} implies
\begin{align*}
  \mathcal A_1(T)+\mathcal A_2(T)\leq CE_0^{2\vartheta}+CE_0^{\frac{\gamma(q_0-2\gamma)}{2(q_0-\gamma)}}\int_0^T\|P-\bar{P}\|_{L^2}^2\mathrm{d}t
  \leq C_5E_0^{2\vartheta}\leq E_0^{\vartheta},
\end{align*}
provided that $E_0\leq\varepsilon_{2,3}:=\min\{\varepsilon_1,C_5^{-1/\vartheta}\}$.
The proof is complete by taking $\varepsilon_2=\min\{\varepsilon_{2,1},\varepsilon_{2,2},\varepsilon_{2,3}\}$.
\end{proof}

We next apply Zlotnik's inequality to improve the $L^{q_0}$ bootstrap bound for the density.

\begin{lemma}\label{l3.4}
Suppose that the bootstrap bounds in \eqref{3.6} hold and $\lambda\geq\lambda_2$, where $\lambda_2>0$ depends only on $\gamma$, $q_0$, and $\mu$.
Then there exists $\varepsilon_3>0$ depending only on $\widehat\rho, \bar\rho, a, \gamma, \mu, \lambda$, and $q_0$ such that
\begin{equation}\label{3.53}
\sup_{0\leq t\leq T}\|\rho\|_{L^{q_0}}\leq\frac{7}{4}\widehat{\rho},
\end{equation}
provided that $E_0\leq \varepsilon_3M^{-6/\vartheta}$.
\end{lemma}

\begin{proof}
The continuity equation \eqref{1.1}$_1$ and the decomposition \eqref{2.7} give
\begin{equation}\label{3.54}
  \frac{\mathrm{d}}{\mathrm{d}t}\int\rho^{q_0}\mathrm{d}\mathbf{x}=-(q_0-1)\int\rho^{q_0}\divv\mathbf{u}\mathrm{d}\mathbf{x}
  =-(q_0-1)\int\rho^{q_0}\divv\mathbf{u}_r\mathrm{d}\mathbf{x}-\frac{q_0-1}{2\mu+\lambda}\int\rho^{q_0}\big(F_s+P-\bar P\big)\mathrm{d}\mathbf{x}.
\end{equation}
Using H\"older's inequality and the estimate for the harmonic component $F_s$ in \eqref{2.10}, we deduce
\begin{align}\label{3.55}
  &\frac{\mathrm{d}}{\mathrm{d}t}\|\rho\|_{L^{q_0}}^{q_0}+\frac{a(q_0-1)}{2\mu+\lambda}
  \left(1-\frac{C_{\gamma,q_0}\mu}{2(3\mu+\lambda)}\right)\|\rho\|_{L^{q_0+\gamma}}^{q_0+\gamma}
  \notag\\ &\leq
  (q_0-1)\|\nabla\mathbf{u}_r\|_{L^\infty}\|\rho\|_{L^{q_0}}^{q_0}+\frac{q_0-1}{2\mu+\lambda}\big(a\bar\rho^\gamma+C_6E_0^{\gamma/2}\big)\|\rho\|_{L^{q_0}}^{q_0}
  \notag\\ &\leq
  (q_0-1)\|\nabla\mathbf{u}_r\|_{L^\infty}\|\rho\|_{L^{q_0}}^{q_0}+\frac{9a(q_0-1)}{8(2\mu+\lambda)}|\Omega|^{-\frac{\gamma}{q_0}}\widehat\rho^\gamma\|\rho\|_{L^{q_0}}^{q_0},
\end{align}
provided that
\[
  E_0\leq\varepsilon_{3,1}:=\min\left\{\varepsilon_2,\left(\frac{a\widehat\rho^\gamma}{8C_6|\Omega|^{\gamma/q_0}}\right)^{2/\gamma}\right\}.
\]
Choose $\lambda_2>0$ such that $C_{\gamma,q_0}\mu/(3\mu+\lambda_2)\le1/4$.
For $\lambda\ge\lambda_2$, define
\[
  y(t):=\|\rho\|_{L^{q_0}}^{q_0}, \quad
  f(y):=\frac{a(q_0-1)}{8(2\mu+\lambda)}|\Omega|^{-\frac{\gamma}{q_0}}\Big(9\widehat\rho^\gamma y-7y^{1+\frac{\gamma}{q_0}}\Big),\quad
  b(t):=(q_0-1)\int_0^t\|\nabla\mathbf{u}_r\|_{L^\infty}y(\tau)\mathrm{d}\tau.
\]
Then H\"older's inequality reduces the estimate \eqref{3.55} to
\begin{equation}\label{3.56}
  y'(t)\leq f\big(y(t)\big)+b'(t).
\end{equation}

To estimate the increments of $b(t)$, we first establish a short-time weighted estimate for $\nabla\dot{\mathbf u}$.
Taking $m=1$ in \eqref{3.46} and integrating over $(0,\sigma(T))$, it follows from \eqref{3.19} and Lemma \ref{l3.2} that
\begin{align*}
  &\sup_{0\leq t\leq \sigma(T)}\big(\sigma\|\sqrt{\rho} \dot{\mathbf{u}}\|_{L^2}^2\big)+\int_0^{\sigma(T)}\sigma\|\nabla \dot{\mathbf{u}}\|_{L^2}^2\mathrm{d}t
  \notag\\ &\leq
  C\int_0^{\sigma(T)}\|\sqrt{\rho}\dot{\mathbf{u}}\|_{L^2}^2\mathrm{d}t+CE_0
  +C\int_0^{\sigma(T)}\sigma\|P-\bar P\|_{L^4}^4\mathrm{d}t
  +C\int_0^{\sigma(T)}\sigma\big(\|\nabla \mathbf{u}_r\|_{L^4}^4+\|\nabla \mathbf{u}_s\|_{L^4}^4\big)\mathrm{d}t
  \notag\\ &\leq
  C\int_0^{\sigma(T)}\|\sqrt{\rho}\dot{\mathbf{u}}\|_{L^2}^2\mathrm{d}t+CE_0+C\sup_{0\leq t\leq T}\|P-\bar P\|_{L^4}^4\int_0^{\sigma(T)}1\mathrm{d}t
  \notag\\ &\quad
  +C\sup_{0\leq t\leq \sigma(T)}\big(\sigma^\frac12\|\nabla\mathbf{u}\|_{L^2}\big) \sup_{0\leq t\leq \sigma(T)}\big(\sigma^\frac12\|\sqrt{\rho} \dot{\mathbf{u}}\|_{L^2}\big)
  \int_0^{\sigma(T)}\|\sqrt{\rho}\dot{\mathbf{u}}\|_{L^2}^2\mathrm{d}t
  \notag\\ &\leq
  CM^2+CE_0+CE_0^{\frac{\gamma(q_0-4\gamma)}{2(q_0-\gamma)}}+CM^2\mathcal A_1^\frac12(T)\sup_{0\leq t\leq \sigma(T)}\big(\sigma^\frac12\|\sqrt{\rho} \dot{\mathbf{u}}\|_{L^2}\big)
  \notag\\ &\leq
  CM^2+CM^4E_0^{\vartheta}+\frac12\sup_{0\leq t\leq \sigma(T)}\big(\sigma\|\sqrt{\rho} \dot{\mathbf{u}}\|_{L^2}^2\big).
\end{align*}
Hence, for $E_0\le M^{-2/\vartheta}$, absorbing the last term gives
\begin{align}\label{3.57}
   \sup_{0\leq t\leq \sigma(T)}\big(\sigma\|\sqrt{\rho} \dot{\mathbf{u}}\|_{L^2}^2\big)+\int_0^{\sigma(T)}\sigma\|\nabla \dot{\mathbf{u}}\|_{L^2}^2\mathrm{d}t\leq CM^2.
\end{align}

For $0\le t_1<t_2\le T$, the bootstrap bounds in \eqref{3.6}, the Gagliardo--Nirenberg inequalities in Lemma~\ref{l2.2},
and the elliptic estimates \eqref{2.10}--\eqref{2.11}, together with Young's inequality, imply
\begin{align}\label{3.58}
  &|b(t_2)-b(t_1)|\leq C\int_{t_1}^{t_2}\|\nabla\mathbf{u}_r\|_{L^\infty}\mathrm{d}t
  \leq C\int_{t_1}^{t_2} \|\nabla\mathbf{u}_r\|_{L^2}^\frac13\|\nabla^2\mathbf{u}_r\|_{L^4}^\frac23 \mathrm{d}t
  \notag\\ &\leq
  C\int_{t_1}^{t_2} \big(\|\nabla\mathbf{u}\|_{L^2}+\|P-\bar P\|_{L^2}\big)^\frac13\Big(\|\rho\dot{\mathbf{u}}\|_{L^4}+\|\rho\|_{L^{4\gamma}}\|\nabla\phi\|_{L^\frac{4\gamma}{\gamma-1}}\Big)^\frac23\mathrm{d}t
  \notag\\ &\leq
  C\int_{t_1}^{t_2} \big(\|\nabla\mathbf{u}\|_{L^2}+\|P-\bar P\|_{L^2}\big)^\frac13\big(\|\sqrt{\rho}\dot{\mathbf{u}}\|_{L^2}+\|\nabla\dot{\mathbf{u}}\|_{L^2}+\|P-\bar P\|_{L^2}\big)^\frac23\mathrm{d}t
  \notag\\ &\leq
  C\sup_{0\leq t\leq T}\bigg(\sigma^\frac16\|\nabla\mathbf{u}\|_{L^2}^\frac13+E_0^{\frac{\gamma(q_0-2\gamma)}{12(q_0-\gamma)}}\bigg)
  \bigg(\int_0^{\sigma(T)}\big(\sigma\|\sqrt{\rho}\dot{\mathbf{u}}\|_{L^2}^2
  +\sigma\|\nabla\dot{\mathbf{u}}\|_{L^2}^2+\|P-\bar P\|_{L^2}^2\big)\mathrm{d}t\bigg)^\frac13
  \bigg(\int_0^{\sigma(T)}\sigma^{-\frac34}\bigg)^{\frac23}
  \notag\\ &\quad
  +C\sup_{0\leq t\leq T}\big(\sigma^\frac12\|\nabla\mathbf{u}\|_{L^2}+1\big)
  \int_{\sigma(T)}^T\big(\sigma\|\sqrt{\rho}\dot{\mathbf{u}}\|_{L^2}^2+\sigma^3\|\nabla\dot{\mathbf{u}}\|_{L^2}^2+\|P-\bar P\|_{L^2}^2\big)\mathrm{d}t
  +N_1(t_2-t_1)
  \notag\\ &\leq
  C\bigg(\mathcal A_1^\frac16(T)+E_0^{\frac{\gamma(q_0-2\gamma)}{12(q_0-\gamma)}}\bigg)\bigg(\mathcal A_1^\frac13(T)+M^\frac23+E_0^{\frac{\gamma(q_0-2\gamma)}{6(q_0-\gamma)}}\bigg)
  +C\Big(\mathcal A_1^\frac12(T)+1\Big)\bigg(\mathcal A_1(T)+\mathcal A_2(T)+E_0^{\frac{\gamma(q_0-2\gamma)}{2(q_0-\gamma)}}\bigg)
  \notag\\[0.2em] &\quad
  +N_1(t_2-t_1)
  \notag\\ &\leq
  CE_0^\frac{\vartheta}{6}M^\frac23+CE_0^{\vartheta}+N_1(t_2-t_1) \leq C_7E_0^\frac{\vartheta}{18}+N_1(t_2-t_1),
\end{align}
provided that $E_0\le M^{-6/\vartheta}$, where Young's inequality is applied with
\[
  N_1=\frac{a(q_0-1)}{8(2\mu+\lambda)}|\Omega|^{-\frac{\gamma}{q_0}}\widehat\rho^{q_0+\gamma}.
\]

Set $N_0=C_7E_0^\frac{\vartheta}{18},\ \xi^*=\left(\frac32\widehat{\rho}\right)^{q_0}$.
Since $\gamma>1$, for every $\xi\ge\xi^*$,
\[
f(\xi)
\leq
\frac{a(q_0-1)}{8(2\mu+\lambda)}
|\Omega|^{-\frac{\gamma}{q_0}}
\widehat\rho^\gamma\xi
\left[9-7\left(\frac32\right)^\gamma\right]\leq
-\frac{3a(q_0-1)}{16(2\mu+\lambda)}
|\Omega|^{-\frac{\gamma}{q_0}}
\widehat\rho^\gamma\xi
\leq-N_1.
\]
Thus, Lemma \ref{l2.4} gives
\[
  y(t)\leq\max\{y(0),\xi^*\}+N_0\leq\left(\frac32\widehat\rho\right)^{q_0}+C_7E_0^{\vartheta/18}
  \leq\left(\frac74\widehat\rho\right)^{q_0},
\]
provided that $E_0\leq\varepsilon_{3,2}:=\min\left\{\varepsilon_{2},\big(4C_7^{-1}\widehat{\rho}^{q_0}\big)^{\frac{18}{\vartheta}}\right\}$.
This completes the proof if $\varepsilon_3=\min\{\varepsilon_{3,1},\varepsilon_{3,2}\}$.
\end{proof}

\begin{proof}[Proof of Proposition \ref{p3.1}]
Choosing $\lambda_*=\max\{\lambda_1,\lambda_2\}$ and $\varepsilon=\varepsilon_3$, then the improved bounds in \eqref{3.7} follow from Lemmas~\ref{l3.2}--\ref{l3.4}.
\end{proof}

\section{Proof of Theorem \ref{t1.1}}\label{sec4}

With the \textit{a priori} estimates established in Section \ref{sec3} at hand, it is a position to show Theorem \ref{t1.1}.

\begin{proof}[Proof of Theorem \ref{t1.1}]
Let $(\rho_0,\mathbf u_0)$ satisfy the assumptions in Theorem~\ref{t1.1}. For $0<\epsilon<1$, let $J_\epsilon$ be a standard nonnegative mollifier and set
\[
  \rho_0^\epsilon:=(1-\epsilon)J_\epsilon*(\rho_0\mathbf 1_\Omega)+\epsilon\bar\rho,\quad \bar\rho^\epsilon:=\frac1{|\Omega|}\int_\Omega\rho_0^\epsilon\mathrm{d}\mathbf{x}.
\]
Let $\mathbf u_0^\epsilon$ and $\phi_0^\epsilon$ solve, respectively,
\[
 \begin{cases}
  \mathbf u_0^\epsilon-\epsilon\Delta\mathbf u_0^\epsilon=\mathbf u_0, &\text{in }\Omega,\\
  \mathbf u_0^\epsilon=\mathbf0, &\text{on }\partial\Omega,
 \end{cases}
  \quad
 \begin{cases}
  \Delta\phi_0^\epsilon=\rho_0^\epsilon-\bar\rho^\epsilon, &\text{in }\Omega,\\
  \partial_{\mathbf n}\phi_0^\epsilon=0, &\text{on }\partial\Omega,
 \end{cases}
\]
with $\int_\Omega\phi_0^\epsilon\mathrm{d}\mathbf{x}=0$.
By standard properties of mollification and elliptic regularity, the approximate initial data satisfy the regularity assumptions \eqref{3.1}--\eqref{3.2}
and the initial bounds in \eqref{1.11}, with $\rho_0^\epsilon\ge\epsilon\bar\rho>0$. Moreover,
\[
  (\rho_0^\epsilon,\mathbf u_0^\epsilon,\phi_0^\epsilon)\to(\rho_0,\mathbf u_0,\phi_0)
  \quad\text{in } L^{q_0}(\Omega)\times H_0^1(\Omega)\times H^2(\Omega)
\]
as $\epsilon\to0$.
In particular, $\bar\rho^\epsilon\to\bar\rho>0$ and $E_0^\epsilon\to E_0$.
For sufficiently small $\epsilon>0$, let $(\rho^\epsilon,\mathbf u^\epsilon,\phi^\epsilon)$ denote the corresponding local strong solution,
with $\bar\rho$ replaced by $\bar\rho^\epsilon$ in the Poisson equation.
The uniform bounds on the approximate initial data ensure that the constants in Section \ref{sec3} can be chosen independently of $\epsilon$.
By reducing the smallness constant in Theorem \ref{t1.1} if necessary, we may also ensure that $E_0^\epsilon$ satisfies the condition of Proposition~\ref{p3.1}.

Since
\[
  \mathcal A_1(0)=\mathcal A_2(0)=0,\quad\mathcal A_3(0)\le M^2,\quad\|\rho_0^\epsilon\|_{L^{q_0}}\le\widehat\rho,
\]
the bootstrap bounds in \eqref{3.6} hold on $[0,T_1]$ for some $T_1>0$ due to continuity. Set
\begin{equation}\label{4.1}
  T_\epsilon^*:=\sup\left\{T>0:\text{ the strong solution exists on }[0,T]\text{ and satisfies }\eqref{3.6}\right\}.
\end{equation}
We claim that
\begin{equation}\label{4.2}
  T_\epsilon^*=\infty.
\end{equation}
Indeed, Proposition \ref{p3.1} yields the improved bounds in \eqref{3.7} for every $T<T_\epsilon^*$, uniformly in $\epsilon$ and $T$.
If $T_\epsilon^*<\infty$, the continuation of the local strong solution under these bounds extends \eqref{3.6} beyond $T_\epsilon^*$,
contradicting \eqref{4.1}. Hence, \eqref{4.2} follows.

For any fixed $0<\tau<T<\infty$, the global approximate solutions satisfy
\begin{equation}\label{4.3}
 \begin{cases}
  \rho^\epsilon\in L^\infty(0,T;L^{q_0}(\Omega)),\quad \phi^\epsilon\in L^\infty(0,T;W^{2,q_0}(\Omega)),\quad
  \sqrt{\rho^\epsilon}\mathbf u^\epsilon\in L^\infty(0,T;L^2(\Omega)),
  \\[0.2em]
  \mathbf u^\epsilon\in L^2(0,T;H_0^1(\Omega))\cap L^\infty(\tau,T;H_0^1(\Omega)),\ \
  \sqrt{\rho^\epsilon}\dot{\mathbf u}^\epsilon\in L^\infty(\tau,T;L^2(\Omega)),\ \ \nabla\dot{\mathbf u}^\epsilon\in L^2(\tau,T;L^2(\Omega)),
 \end{cases}
\end{equation}
uniformly in $\epsilon$. Thus, up to a subsequence,
\[
  \rho^\epsilon\to\rho\quad\text{weakly* in }L^\infty(0,T;L^{q_0}(\Omega)),\quad
  \mathbf u^\epsilon\to\mathbf u \quad\text{weakly in }L^2(0,T;H_0^1(\Omega)).
\]
Combining the standard compactness argument (cf. \cite{PL98,NS04,S2020}) with the $L^{q_0}$ bound and Vitali's theorem, we obtain
\[
  P(\rho^\epsilon)\to P(\rho)\quad\text{strongly in }L^q(\Omega\times(0,T)), \quad 1\le q<q_0/\gamma.
\]
Moreover, owing to $\bar\rho^\epsilon\to\bar\rho$, the Neumann elliptic estimates yield
\[
  \phi^\epsilon\to\phi\quad\text{strongly in }L^q(0,T;W^{2,q}(\Omega)),\quad 1<q<q_0.
\]

The estimates in \eqref{4.3} and the preceding convergences allow us to pass to the limit in the weak formulations by standard compactness arguments.
We thus obtain a global weak solution $(\rho,\mathbf u,\phi)$ in the sense of Definition \ref{d1.1} satisfying the estimates stated in Theorem \ref{t1.1}.
\end{proof}

\section*{Conflict of interest}
The authors declare that they have no conflict of interest.

\section*{Data availability}
No datasets were generated or analyzed during the current study.

\end{document}